\documentclass[11pt]{amsart}
\usepackage{mathrsfs}
\usepackage{amssymb}
\usepackage{amsthm}
\usepackage{titletoc}
\usepackage{amssymb}
\usepackage{amsxtra}
\usepackage{mathrsfs}
\usepackage{hyperref}

\newtheorem{theorem}{Theorem}
\newtheorem{corollary}{Corollary}

\newtheorem{lemma}{Lemma}
\newtheorem{remark}{Remark}
\newtheorem{definition}{Definition}
\newtheorem{question}{Question}
\newtheorem{example}{Example}

\title[Meromorphic functions sharing four values]{On the solution of a question for meromorphic functions sharing three values IM and a fourth value CM in the complex plane}
	\author{Xiao-Min Li }
	\address{Department of Mathematics, Ocean University of China, Qingdao, Shandong 266100, P. R. China}
	\email{lixiaomin@ouc.edu.cn(X.-M. Li)}
	\author{Jun Wang}
\address{School of Mathematical Sciences, Fudan University, Shanghai, 200433, P. R. China}
\email{majwang@fudan.edu.cn(J. Wang)}
\author{Hong-Xun Yi}
\address{Department of Mathematics, Shandong University, Jinan, Shandong 250199, P. R. China}
	\email{hxyi@sdu.edu.cn(H.-X. Yi)}
		\thanks{Project supported in part by the NSF of Shandong Province, China (No. ZR2019MA029),
		the FRFCU(No. 3016000842464005).}
\begin{document}
\begin{abstract}
In 1926, Nevanlinna (Acta Math. {\bf 48} (1926), no. 3-4, 367-391) proved that if two
distinct and non-constant meromorphic functions $f$ and $g$ share $a_1,$ $a_2,$ $a_3,$ $a_4$ CM in the complex plane, where $a_1,$ $a_2,$ $a_3,$ $a_4$ are four distinct complex values in the extended complex plane, then $f$ is a M\"{o}bius transformation of $g$, two of the shared values, say $a_1$ and $a_2$, are Picard exceptional values, and the cross ratio $\left(a_1, a_2, a_3, a_4\right)=-1.$ Later on,
Gundersen (Trans. Amer. Math. Soc. {\bf 277} (1983), no. 2, 545-567. Correction: {\bf 304} (1987), 847-850) proved that if two non-constant meromorphic functions share two values CM and share two other values IM, then the functions share all four values CM. In 1992, Gundersen (Complex Var. Elliptic Equ. {\bf 20} (1992), no. 1-4, 99-106) proposed the following famous question: if two non-constant meromorphic functions share three distinct values IM and share a fourth value CM, then do the functions necessarily share all four values CM? This question is a long-standing question in the studies of the Nevanlinna's value distribution theory of meromorphic functions, and has not been completely resolved by now. In this paper, we completely resolve the question, and prove that if two distinct and non-constant meromorphic functions $f$ and $g$ share $a_1,$ $a_2,$ $a_3$ IM and $\infty$ CM, where $a_1,$ $a_2,$ $a_3$ are three distinct and finite complex values in the complex plane, then $f$ and $g$ also share $a_1,$ $a_2,$ $a_3$ CM. The main result in this paper also improves the corresponding results from Mues (Complex Var. Elliptic Equ. {\bf 12} (1989), no. 1-4, 169-179), Gundersen (Complex Var. Elliptic Equ. {\bf 20} (1992), no. 1-4, 99-106) and Steinmetz (Southeast Asian Bull.Math. {\bf 36} (2012), no. 3, 399-417). A famous example provided in (J. London Math. Soc. {\bf 20} (1979), no. 2, 457-466) shows that the results obtained in this paper, in a sense, are best possible.
\end{abstract}
\vskip 2mm
\par
\maketitle
\section{Introduction and main results}
\vskip 2mm
\par In 1926, Nevanlinna \cite{Nevanlinna1926} proved the following theorems that are called the five-value theorem and the four-value theorem respectively:
\begin{theorem}\rm{}(\cite[p.109]{Nevanlinna1926})\label{TheoremA}
	If $f$ and $g$ are non-constant meromorphic functions that share five values IM, then $f=g.$
\end{theorem}
\begin{theorem}\rm{}(\cite[p.122]{Nevanlinna1926})\label{TheoremB}
	Let $f$ and $g$ be distinct non-constant meromorphic functions, and let $a_1,$ $a_2,$ $a_3,$ $a_4$ be four distinct complex values in the extended complex plane. If $f$ and $g$ share $a_1,$ $a_2,$ $a_3,$ $a_4$ CM, then $f$ is a M\"{o}bius transformation of $g$, two of the shared values, say $a_1$ and $a_2$, are Picard exceptional values, and the cross ratio $\left(a_1, a_2, a_3, a_4\right)=-1.$
\end{theorem}
\par\vskip 2mm In 1979, Gundersen \cite[p.457]{Gundersen1979} wrote: L. Rubel has asked what can be said if CM is replaced by IM in the hypothesis of
Theorem \ref{TheoremB}?  In this direction, Gundersen \cite{Gundersen1979} and Gundersen \cite{Gundersen1983} proved the following results respectively, which improved Theorem \ref{TheoremB}:
\begin{theorem}\rm{}(\cite[Theorem 1]{Gundersen1979}) \label{TheoremC}
	If $f$ and $g$ share three values CM and share a fourth value IM, then they share all four values CM.
\end{theorem}
\begin{theorem}\rm{}(\cite[Theorem 1]{Gundersen1983})\label{TheoremD}
	If two non-constant meromorphic functions share two values CM and share two other values IM, then the functions share all four values CM.
\end{theorem}
\par
\vskip 2mm
The following famous example given in Gundersen \cite{Gundersen1979} shows that the Rubel's question above is negative:
\begin{example}\rm{}(\cite[pp.458-459]{Gundersen1979}).\label{Example1.1}
	Let $f=\frac{e^{h}+1}{(e^{h}-1)^2}$ and $g=\frac{(e^{h}+1)^2}{8(e^{h}-1)},$ where $h$ is a non-constant entire function, we can easily verify that $f$ and $g$ share $0,$ $1,$ $\infty,$ $-1/8$ DM, and $f$ is not a M\"{o}bius transformation of $g.$
\end{example}
\vskip 2mm
\par Based upon Theorem \ref{TheoremD} and Example \ref{Example1.1}, Gundersen \cite{Gundersen1992} proposed the following famous question:
\begin{question}\rm{}(\cite[p.100]{Gundersen1992}).\label{Question1.1}
	If two non-constant meromorphic functions share three values IM and share a fourth value CM, then do the functions necessarily share all the four values CM?
\end{question}
\par
\vskip 2mm 
It is well-known that Question \ref{Question1.1} is the main open question in the uniqueness theory of meromorphic functions. This open question is a long-standing and difficult question, and has not been completely resolved by now. As far as we know, there are many other research works concerning Question 1.7, such as Mues \cite{Mues1989F1,Mues1989F2}, Gundersen \cite{Gundersen1992}, Reinders \cite{Reinders1992}, Ueda \cite{Ueda1994}, Wang \cite{Wang2001}, Steinmetz \cite{Steinmetz1988,Steinmetz2012}, Huang \cite{Huang2004} and Huang-Du \cite{Huang2004D2}, etc. In addition, Ishizaki \cite{Ishizaki2000}, Li \cite{Li2001} and Yao \cite{Yao2003} studied the uniqueness question of non-constant meromorphic functions sharing four distinct small functions.
\vskip 2mm
\par
For the ease of the following narrations, we introduce the following notation: let $f$ and $g$ be two non-constant meromorphic functions, and let $a$ be a value in the extended complex plane. Suppose that $f$ and $g$ share $a$ IM. Next we denote by $\overline{N}_{E}(r,a;f,g)$ the reduced counting function of those common $a$-points of $f$ and $g$ in $\left|z\right|<r,$ where each point in $\overline{N}_{E}(r,a;f,g)$ has the same multiplicities regard to $f$ and $g,$ and each point in $\overline{N}_{E}(r,a;f,g)$ is counted only once. We define
$\tau(a)=\liminf\limits_{r\rightarrow\infty}\frac{\overline{N}_{E}(r,a;f,g)}{\overline{N}(r,a,f)},$ if $\overline{N}(r,a,f)\not=0$
for the large positive number $r,$ where $\overline{N}(r,a,f)=\overline{N}(r,f)$ if $a=\infty.$ From the definition we see that $0\le \tau(a)\le 1.$ We mention that $\tau(a)$ is defined as $\tau(a)=1,$ when $f$ and $g$ share $a$ CM. We recall the following result due to Mues \cite{Mues1989F2} that improved Theorem \ref{TheoremD}:
\begin{theorem}\rm{}(\cite[Theorem 1]{Mues1989F2}). \label{TheoremE}
	Let $f$ and $g$ be two distinct non-constant meromorphic functions in the complex plane that share four distinct values in the extended complex plane. Suppose that one of the four distinct values is shared by $f$ and $g$ CM and that $\tau(b)>2/3$ for another value $b$ of the four distinct shared values. Then $f$ and $g$ share all the four distinct values CM.
\end{theorem}
\par
\vskip 2mm
	In 1992, Gundersen\cite{Gundersen1992} studied Question \ref{Question1.1}, and proved the following result:
\begin{theorem}\rm{}(\cite[Theorem 1]{Gundersen1992}). \label{TheoremF}
	Let $f$ and $g$ be two distinct non-constant meromorphic functions that share $a_1,$ $a_2,$ $a_3$ IM and $\infty$ CM, where $a_1,$ $a_2,$ $a_3$ are three distinct and finite complex values in the complex plane. Suppose that there exist some real constant $\sigma>4/5,$ and some set $I\subset (0,+\infty)$ that has infinite linear measure such that $N(r,f) \ge \sigma T(r,f)$ for $r\in I.$ Then $f$ and $g$ also share $a_1,$ $a_2,$ $a_3$ CM.
\end{theorem}
\par
\vskip 2mm Later on, Steinmetz\cite{Steinmetz2012} also studied Question \ref{Question1.1}, and got the following result:
\par
\vskip 2mm
\begin{theorem}\rm{}(\cite[pp.411-412, key observations]{Steinmetz2012}). \label{TheoremG}
	Let $f$ and $g$ be two distinct transcendental meromorphic functions that share $a_1,$ $a_2,$ $a_3,$ $\infty$ IM, where $a_1,$ $a_2,$ $a_3$ are three distinct and finite complex values in the complex plane, and let
\begin{equation}\label{eq1.1X1}
\Phi_{f}=\frac{f'(f-g)}{(f-a_1)(f-a_2)(f-a_3)} \, \, \text{and} \, \, \Phi_{g}=\frac{g'(f-g)}{(g-a_1)(g-a_2)(g-a_3)}.
\end{equation}
Then we have the following conclusions:
\par
\vskip 2mm  (i) If additionally $f$ and $g$ share $\infty$ CM, then $\Phi_{f}$ and $\Phi_{g}$ are entire functions such that
$N\left(r, 1/\Phi_{f}\right)+N\left(r,1/\Phi_{g}\right)=S(r,f).$
\par
\vskip 2mm  (ii) If additionally $f$ and $g$ are distinct and non-constant meromorphic functions of finite lower order, then
 $\Phi_{f}=P_fe^{Q}$ and $\Phi_{g}=P_ge^{-Q},$ where $P_f$ and $P_g$ are non-vanishing polynomials, and $Q$ is a non-constant polynomial such that
\begin{equation}\nonumber
\Psi=\Phi_{f}\Phi_{g}=\frac{f'g'{(f-g)}^2}{(f-a_1)(f-a_2)(f-a_3)(g-a_1)(g-a_2)(g-a_3)}
\end{equation}
and
\begin{equation}\nonumber
T(r,f)=T(r,g)+O(\log r)=r^{\deg(Q)}\left(1+o(1)\right),
\end{equation}
as $r\rightarrow+\infty.$ Here and in what follows, $\Psi$ is an entire function as defined in Lemma \ref{Gundersen1992} in Section 2 of this paper such that
$T(r,\Psi)=O(\log r),$ as $r\rightarrow+\infty.$
\end{theorem}
\vskip 2mm
\par In this paper, we will prove the following result that improves  Theorem \ref{TheoremB}-Theorem \ref{TheoremG} and
gives an affirmative answer to Question \ref{Question1.1} completely:
\vskip 2mm
\par
\begin{theorem}\rm{}\label{Theorem1.2}
Let $f$ and $g$ be two distinct non-constant meromorphic functions, and let $a_1,$ $a_2,$ $a_3$ be three distinct and finite values in the complex plane. Suppose that $f$ and $g$ share $a_1,$ $a_2,$ $a_3$ IM and $\infty$ CM. Then $f$ and $g$ also share $a_1,$ $a_2,$ $a_3$ CM.
\end{theorem}
\vskip 2mm
\par From Theorem \ref{Theorem1.2} we get the following result:
\vskip 2mm
\par
\begin{corollary}\rm{}\label{Corollary1.1}
Let $f$ and $g$ be two distinct non-constant entire functions, and let $a_1,$ $a_2,$ $a_3$ be three distinct and finite values in the complex plane. Suppose that $f$ and $g$ share $a_1,$ $a_2,$ $a_3$ IM. Then $f$ and $g$ share $a_1,$ $a_2,$ $a_3$ CM.
\end{corollary}
\vskip 2mm
\par Throughout this paper, by meromorphic functions we will always mean meromorphic functions in the complex plane. We adopt the standard notations in the Nevanlinna's theory of meromorphic functions as explained in \cite{Hayman1964,Laine1993,Yang2003}. We give examples as follows: for a meromorphic function $f$ we denote by  $m(r,f),$ $N(r,f)$ and $T(r,f)$ the proximity function on $|z|=r,$ the counting function of poles of $f$ in $|z|<r$ counted with multiplicities,  and the characteristic function of $f$ respectively, they are defined as
\begin{equation}\nonumber
m(r,f)=\frac{1}{2\pi}\int_{0}^{2\pi}\log^{+}|f(re^{i\theta})|d\theta, \, \, N(r,f)=\int_{0}^{r}\frac{n(t,f)-n(0,f)}{t}dt+n(0,f)\log r
\end{equation}
and $T(r,f)=m(r,f)+N(r,f)$ for $0<r<+\infty$ respectively, where and in what follows, $\log ^{+}x =\log x,$ if $x\geq 1,$ and $\log ^{+}x =0,$ if $0\leq x<1,$ while $n(t,f)$ denotes the number of poles of $f$ in $|z|<t$ counted with multiplicities for the positive constant $t$ satisfying $0<t\leq r,$ and $n(0,f)$ denotes the multiplicities of the origin point of a pole of $f.$
Throughout this paper, we denote by $E\subset[0,+\infty)$ a set of finite linear measure. For a meromorphic function $f$ in the complex plane, we denote by $S(r,f)$ any quantity satisfying $S(r,f)=o(T(r,f)),$ as $r \not\in E$ and $r\rightarrow\infty.$ Following \cite[p.466]{Halburd2006}, we define the following notion and notation: let $f$ be a non-constant meromorphic function. Next we denote by $S(f)$ the set of all meromorphic functions $\alpha$ in the complex plane such that $T(r,\alpha)=o(T(r,f))$ for all $r\in[0,+\infty)$ possibly outside a set $E\subset[0,+\infty)$ of finite linear measure. Functions in the set $S(f)$ are called small compared to $f$, or slowly moving with respect to $f,$ or a small function of $f$ for short.
\vskip 2mm
\par Let $f$ and $g$ be two non-constant meromorphic functions in the complex plane, and let $a$ be a value in the extended plane. We say that $f$ and $g$ share the value $a$ CM, provided that $f$ and $g$ have the same set of $a$-points in the complex plane, and each common $a$-points of $f$ and $g$ have the same multiplicities. We say that $f$ and $g$ share the value $a$ IM, provided that $f$ and $g$ have the same set of $a$-points in the complex plane, and each common $a$-point is counted only once. We say that $f$ and $g$ share the value $a$ DM, provided that $f$ and $g$ share $a$ IM, and at least one common $a$-point of $f$ and $g$ has different multiplicities related to $f$ and $g$ (cf.\cite{Gundersen1979,Gundersen1983}). Here we mention that Mues \cite{Mues1989F1} and Yang-Yi \cite{Yang2003} defined the notion of the DM shared value as follows: we say that $f$ and $g$ share the value $a$ DM, provided that $f$ and $g$ share $a$ IM, and that $f-a$ and $g-a$ have no common zeros of the same multiplicities. For a non-constant meromorphic function $f,$ we denote by $\mu(f)$ and $\rho(f)$ the lower order and the order of $f$ respectively, the relevant definitions can be found, for example, in \cite{Hayman1964,Laine1993, Yang2003}. For convenience, we recall them below:
\begin{definition}\rm{} \label{Definition 1.1} For a non-constant meromorphic function $f,$ the lower order $\mu(f)$ and the order $\rho(f)$ of $f,$ are defined as
	\begin{equation}\nonumber
		\mu(f)=\liminf\limits_{r\rightarrow\infty}\frac{\log {T(r,f)}}{\log r} \quad \text{and} \quad \rho(f)=\limsup\limits_{r\rightarrow\infty}\frac{\log {T(r,f)}}{\log r} \; \text{respectively}.
	\end{equation}
\end{definition}
\par
\vskip 2mm Following \cite[p.293]{Edrei1959}, we say that a meromorphic function $f$ in the complex plane is of regular growth, provided $\mu(f)=\rho(f).$ Here we mention that both of $\mu(f)$ and $\rho(f)$ may be $+\infty.$
\vskip 2mm
\par
\maketitle
\section{Preliminaries}
\vskip 2mm
\par In this section, we will introduce some results that play an important role in proving the main result in this paper. First of all, we introduce the following result due to Adams-Straus\cite{Adams1971}:
\par
\vskip 2mm
\begin{lemma}\rm{}(\cite[Theorem 3]{Adams1971})\label{Adams1971}
	Let $f$ and $g$ be two non-constant rational functions, and let $a_1,$ $a_2,$ $a_3,$ $a_4$ be four distinct values in the extended complex plane. If $f$ and $g$ share $a_1,$ $a_2,$ $a_3,$ $a_4$ IM, then $f=g.$
\end{lemma}
\vskip 2mm
\par In the early 1980s, E. Mues discovered an auxiliary function named Mues$'$s function that plays an important role in the studies of the value distribution theory of meromorphic functions sharing four values by now. We recall the following result due to Gundersen \cite{Gundersen1992}:
\vskip 2mm
\par
\begin{lemma}\rm{}(\cite[Lemma 2]{Gundersen1992} and \cite[pp.176-177]{FrankHua1999})\label{Gundersen1992}
	Let $f$ and $g$ be two distinct non-constant meromorphic functions, and let $a_1,$ $a_2,$ $a_3$ be three distinct and finite values in the complex plane. If $f$ and $g$ share $a_1,$ $a_2,$ $a_3,$ $\infty$ IM, then the function
	\begin{equation}\label{eq2.1X19}
		\Psi=\frac{f'g'{(f-g)}^2}{(f-a_1)(f-a_2)(f-a_3)(g-a_1)(g-a_2)(g-a_3)}
	\end{equation}
is an entire function such that its Nevanlinna's characteristic function satisfies
	\begin{equation}\nonumber
		T(r,\Psi)=S(r,f).
	\end{equation}
\end{lemma}
\par
\vskip 2mm
The following result was proved by Nevanlinna \cite{Nevanlinna1926} originally:
\par
\vskip 2mm
\begin{lemma}\rm{}(\cite[p.373, Satz 3]{Nevanlinna1926} or \cite[Theorem 2]{Gundersen1979})\label{Nevanlinna1926}
		Let $f$ and $g$ be two distinct non-constant meromorphic functions, and let $a_1,$ $a_2,$ $a_3$ be three distinct and finite values,  $f$ and $g$ share $a_1,$ $a_2,$ $a_3$ and $\infty$ IM. Then for $r>0$ we have the following assertions:
\par
\vskip 2mm
\noindent{}
(i) $T(r,f)=T(r,g)+S(r,f)+S(r,g);$
\par
\vskip 2mm\noindent{}(ii) $2T(r,f)=\overline{N}(r,f)+\sum\limits_{j=1}^{3}\overline{N}\left(r,\frac{1}{f-a_{j}}\right)+S(r,f)+S(r,g);$
\par
\vskip 2mm\noindent{}(iii) $2T(r,g)=\overline{N}(r,g)+\sum\limits_{j=1}^{3}\overline{N}\left(r,\frac{1}{g-a_{j}}\right)+S(r,f)+S(r,g).$
	\end{lemma}
\par
\vskip 2mm The following result and its proof can be found in Gundersen \cite{Gundersen1983}:
\vskip 2mm
\par
\begin{lemma}\rm{}(\cite[p.550, Lemma 2]{Gundersen1983} and \cite[p.549, Corollary 1(II)]{Gundersen1983}) \label{Gundersen1983}
Let $f$ and $g$ be two distinct non-constant meromorphic functions that share four distinct values $a_1,$ $a_2,$ $a_3$ and $a_4$ IM, where $a_4=\infty.$ Then the following statement holds:
\vskip 2mm
\par (i) $N_0(r,0,f')+N_0(r,0,g')=S(r,f)+S(r,g),$ where $N_0(r,0,f')$ and $N_0(r,0,g')$ ``count'' respectively only those points in $N_0(r,0,f')$ and $N_0(r,0,g')$ which do not occur when $f(z)=g(z)=a_j$ for some $j\in\{1,2,3,4\}.$
\vskip 2mm
\par (ii) $\sum\limits_{j=1}^{4}N^*(r,a_j;f,g)=S(r,f)+S(r,g),$ where $N^*(r,a_j;f,g)$ with $j\in\{1,2,3,4\}$ denotes the counting function of those common $a_j$-points of $f$ and $g$ in $\left|z\right|<r$ that are multiple for both $f$ and $g$, and each such point in $N^*(r,a_j;f,g)$ is counted according to the number of the times of the smaller of the two multiplicities.
\end{lemma}
\vskip 2mm
\par
The following result is due to Li-Yi\cite{LiYi2007}:
\vskip 2mm
\par
\begin{lemma}\rm{}(\cite[Lemma 9]{LiYi2007})\label{LiYi2007}
Let $f$ and $g$ be two distinct and non-constant meromorphic
functions. Suppose that $f$ and $g$ share $a_1,$ $a_2,$ $a_3,$ $a_4$ IM,
where $a_1,$ $a_2,$ $a_3$ are three distinct finite values and $a_4=\infty.$ Then we have
\begin{equation}\label{eq2.2X19}
\frac {1}{7} T(r,f)\leq\sum\limits_{j=1}^{3}\sum\limits_{l=1}^{6}\overline{N}_{(1,l)}\left(r,a_j;f,g\right)
+\sum\limits_{j=1}^{3}\sum\limits_{m=2}^{6}\overline{N}_{(m,1)}\left(r,a_j;f,g\right)
+S(r,f)+S(r,g)
\end{equation}
and
\begin{equation}\label{eq2.3X19}
\frac {1}{7} T(r,g)\leq\sum\limits_{j=1}^{3}\sum\limits_{l=1}^{6}\overline{N}_{(1,l)}\left(r,a_j;g,f\right)
+\sum\limits_{j=1}^{3}\sum\limits_{m=2}^{6}\overline{N}_{(m,1)}\left(r,a_j;g,f\right)
+S(r,f)+S(r,g).
\end{equation}
Here $\overline{N}_{(1,l)}\left (r,a;f,g \right )$ denotes the reduced counting function of those common zeros of $f-a$ and $g-a$ in $|z|<r$ that are simple zeros of $f-a$ and zeros of $g-a$ of multiplicity $l$ for $a\in\{a_1,a_2,a_3,\infty\},$ while $\overline{N}_{(m,1)}\left(r,a;f,g\right)$ denotes the reduced counting function of those common zeros of $f-a$ and $g-a$ in $|z|<r$ that are zeros of $f-a$ of multiplicity $m$ and simple zeros of $g-a$ for $a\in\{a_1,a_2,a_3,\infty\}.$
 \end{lemma}
\par
\vskip 2mm
The following result proved by Heittokangas-Korhonen-R\"{a}tty\"{a} \cite{Heittokangas2004} is a generalization of the estimate of the proximity function of the logarithmic derivative of a non-vanishing meromorphic function due to Gol$'$dberg-Grinshtein\cite{Goldberg1976}, and a generalization of the Nevanlinna lemma on the logarithmic derivative due to Hiong \cite{Hiong1955}:
\par
\vskip 2mm
\begin{lemma}\rm{}(\cite[Theorem 4.1]{Heittokangas2004})\label{Heittokangas2004} Let $k$ and $j$ be integers such that $k > j\geq 0,$ and let $f$ be a meromorphic function in the complex plane such that $f^{(j)}$ does not vanish identically. Then there exists a positive constant $r_0 > 1$ such that
\begin{equation}\nonumber
m\left(r,\frac{f^{(k)}}{f^{(j)}}\right)\leq (k-j)\log^{+}\left(\frac{T(\rho,f)}{r}\cdot\frac{\rho}{\rho-r}\right) +\log\frac{k!}{j!}+(k-j)5.3078
\end{equation}
for all $r_0<r<\rho$
\end{lemma}
\vskip 2mm
\par
\begin{remark}\rm{}\label{Remark2.1} A proof of Lemma \ref{Nevanlinna1926} can be found in \cite{Gundersen1992}, where \cite[Theorem 2.1]{Hayman1964} was used. From \cite[Theorem 2.1]{Hayman1964} and the proof of \cite[Lemma 1]{Gundersen1992} we see that the error term $S(r,f)+S(r,g)$ in Lemma \ref{Nevanlinna1926} can be expressed as
\begin{equation}\label{eq2.4X19}
\begin{split}
 S(r,f)+S(r,g)&=m\left(r,\frac{f'}{f}\right)+m\left(r,\frac{g'}{g}\right)+m\left(r,\sum\limits_{j=1}^3\frac{f'}{f-a_j}\right)\\
 \end{split}
 \end{equation}
 \begin{equation}\nonumber
 \begin{split}
 &\quad +m\left(r,\sum\limits_{j=1}^3\frac{g'}{g-a_j}\right)+O(1) \, \, \text{for} \, \,  r>0.
\end{split}
\end{equation}
\end{remark}
\vskip 2mm
\par
\begin{remark}\rm{}\label{Remark2.2}
A proof of Lemma \ref{Gundersen1983} can be found in \cite{Gundersen1992}, where the Mues$'$s function $\Psi$ of \eqref{eq2.1X19} was used. From the proof of \cite[Lemma 3]{Gundersen1992} and the decomposition
\begin{equation}\label{eq2.5X19}
\Psi=\sum\limits_{s=1}^3\sum\limits_{t=1}^3\frac{c_{st}f'}{f-a_s}\cdot\frac{g'}{g-a_t}
\end{equation}
 (cf.\cite[pp.176-177]{FrankHua1999}), we see that the error term $S(r,f)+S(r,g)$ in Lemma \ref{Gundersen1983} can be expressed as
\begin{equation}\nonumber 
S(r,f)+S(r,g)=m(r,\Psi)+O(1)= m\left(r,\sum\limits_{s=1}^3\sum\limits_{t=1}^3\frac{c_{st}f'}{f-a_s}\cdot\frac{g'}{g-a_t}\right)+O(1) \, \, \text{for} \, \, r>0.
\end{equation}
Here and in what follows, $c_{st}$ is some complex constant for any pair of positive integers $s$ and $t$ satisfying $1\leq s\leq 3$ and  $1\leq t\leq 3.$
\end{remark}
\vskip 2mm
\par
\begin{remark} \rm{}\label{remark2.3} From the context of \cite{LiYi2007} we see that Lemma \ref{Nevanlinna1926} and Lemma \ref{Gundersen1983} were used for proving \cite[Lemma 9]{LiYi2007}. Therefore, it follows from  Remarks \ref{Remark2.1} and \ref{Remark2.2} that the error term $S(r,f)+S(r,g)$ in \eqref{eq2.2X19} and \eqref{eq2.3X19} of Lemma \ref{LiYi2007} can be expressed as
\begin{equation}\label{eq2.6X19}
\begin{split}
&\quad S(r,f)+S(r,g)= m\left(r,\frac{f'}{f}\right)+m\left(r,\sum\limits_{j=1}^3\frac{f'}{f-a_j}\right)+m\left(r,\frac{g'}{g}\right)\\
&\quad +m\left(r,\sum\limits_{j=1}^3\frac{g'}{g-a_j}\right)+m\left(r,\sum\limits_{s=1}^3\sum\limits_{t=1}^3\frac{c_{st}f'}{f-a_s}\cdot\frac{g'}{g-a_t}\right)+O(1) \, \, \text{for} \, \, r>0.
\end{split}
\end{equation}
\end{remark}
\vskip 2mm
\par
\begin{remark} \rm{}\label{remark2.4} Suppose that $f'\not\equiv 0$ and $g'\not\equiv 0.$  Then, from Lemma \ref{Heittokangas2004} and \eqref{eq2.6X19} in Remark \ref{remark2.3} we deduce that the error term $S(r,f)+S(r,g)$ of \eqref{eq2.2X19} and \eqref{eq2.3X19} in Lemma \ref{LiYi2007} can be estimated as
\begin{equation}\label{eq2.7X19}
\begin{split}
&\quad S(r,f)+S(r,g)\leq m\left(r,\frac{f'}{f}\right)+m\left(r,\sum\limits_{j=1}^3\frac{f'}{f-a_j}\right)+m\left(r,\frac{f''}{f'}\right) +m\left(r,\frac{g'}{g}\right)\\
&\quad+m\left(r,\sum\limits_{j=1}^3\frac{g'}{g-a_j}\right)+m\left(r,\frac{g''}{g'}\right)+m\left(r,\sum\limits_{s=1}^3\sum\limits_{t=1}^3\frac{c_{st}f'}{f-a_s}
\cdot\frac{g'}{g-a_t}\right)+O(1)\\
&\leq m\left(r,\frac{f'}{f}\right)+m\left(r,\frac{f''}{f'}\right)+m\left(r,\frac{g'}{g}\right)+m\left(r,\frac{g''}{g'}\right)\\
 &\quad +4\sum\limits_{s=1}^3\left(m\left(r,\frac{f'}{f-a_s}\right)+m\left(r,\frac{g'}{g-a_s}\right)\right)+O(1)\\
&\leq C_1\log^{+}\left(\frac{T(\rho,f)}{r}\cdot\frac{\rho}{\rho-r}\right)+C_1\log^{+}\left(\frac{T(\rho,g)}{r}\cdot\frac{\rho}{\rho-r}\right)+C_1\\
\end{split}
\end{equation}
for any pair of positive numbers $r$ and $\rho$ such that $1<r_0<r<\rho<R.$ Here $C_1$ and $r_0$ are some two large positive constants such that
$\min\{C_1,r_0\}>1.$ On the other hand, from \eqref{eq2.7X19}, Lemma \ref{Nevanlinna1926} (i) and the expression \eqref{eq2.4X19} of the error term $S(r,f)+S(r,g)$ in Remark \ref{Remark2.1} we have
\begin{equation}\label{eq2.8X19}
\begin{split}
T(\rho,g)&\leq T(\rho,f)+C_1\left( \log^{+}\left( \frac{T(\rho,f)}{\rho}\cdot\frac{2\rho}{2\rho-\rho}\right)
+\log^{+}\left( \frac{T(\rho,g)}{\rho}\cdot\frac{2\rho}{2\rho-\rho}\right)+1\right)\\
&\leq T(\rho,f)+C_2\left( \log^{+}T(\rho,f)+\log^{+}T(\rho,g)+1\right)\\
\end{split}
\end{equation}
for the pair of positive numbers $r$ and $\rho$ such that $1<r_0<r<\rho<2\rho<R<+\infty.$ Noting that $f$ and $g$ are a non-constant meromorphic function, we have that $T(r,f)$ and $T(r,g)$ are two monotonically increasing, continuous and unbounded functions for $0 < r < +\infty$ (cf.\cite[p.8, Corollary 1]{Hayman1964}), we see that there exists some positive constant, say $r_0>1$ such that
\begin{equation}\label{eq2.9X19}
\log^{+}T(\rho,f)=\log T(\rho,f) \, \, \text{and} \, \,  C_2\log^{+}T(\rho,g)+C_2=C_2\log T(\rho,g)+C_2\leq \frac{1}{2}T(\rho, g)
  \end{equation}
 for the pair of positive numbers $r$ and $\rho$ such that $1<r_0<r<\rho<2\rho<R<+\infty$ in \eqref{eq2.8X19}. From \eqref{eq2.8X19} and \eqref{eq2.9X19} we have
\begin{equation}\label{eq2.10X19}
T(\rho, g)\leq 2T(\rho,f)+2C_2\log^{+}T(\rho, f)=2T(\rho,f)+2C_2\log T(\rho,f)
\end{equation}
for the pair of positive numbers $r$ and $\rho$ such that $1<r_0<r<\rho<2\rho<R<+\infty$ in \eqref{eq2.8X19}. We substitute \eqref{eq2.10X19} into \eqref{eq2.7X19}, and then we have
\begin{equation}\label{eq2.11X19}
\begin{split}
&\quad S(r,f)+S(r,g)\leq m\left(r,\frac{f'}{f}\right)+ m\left(r,\frac{f''}{f'}\right)+m\left(r,\sum\limits_{j=1}^3\frac{f'}{f-a_j}\right)+ m\left(r,\frac{g'}{g}\right)\\
&\quad +m\left(r,\frac{g''}{g'}\right)+m\left(r,\sum\limits_{j=1}^3\frac{g'}{g-a_j}\right)+m\left(r,\sum\limits_{s=1}^3\sum\limits_{t=1}^3\frac{c_{st}f'}{f-a_s}\cdot\frac{g'}{g-a_t}\right)+O(1)\\
&\leq C_1\log^{+}\left(\frac{T(\rho,f)}{r}\cdot\frac{\rho}{\rho-r}\right)+C_1\log^{+}T(\rho,g)
+C_1\log^{+}\frac{\rho}{r(\rho-r)}+C_1\\
&\leq C_1\log^{+}\left(\frac{T(\rho,f)}{r}\cdot\frac{\rho}{\rho-r}\right)+C_1\log^{+}(2T(\rho,f)+2C_2\log T(\rho,f))\\
&\quad +C_1\log^{+}\frac{\rho}{r(\rho-r)}+C_1\leq C_1\left(\log^{+}T(\rho,f)+\log^{+}\frac{1}{r}+\log^{+}\frac{\rho}{\rho-r}\right)\\
&\quad+C_1\left(\log^{+}T(\rho,f)+\log^{+}\log^{+}T(\rho,f)+\log^{+}C_2+3\log 2\right)\\
&\quad +C_1\left(\log^{+}\frac{1}{r}+\log^{+}\frac{\rho}{\rho-r}\right)=2C_1\log^{+}T(\rho,f)+C_1\log^{+}\log^{+}T(\rho,f)\\
&\quad+2C_1\left(\log^{+}\frac{1}{r}+\log^{+}\frac{\rho}{\rho-r}\right)+C_1\left(\log^{+}C_2+3\log 2\right)\leq 2C_1\log^{+}T(\rho,f)\\
&\quad +C_1\log^{+}\log^{+}T(\rho,f)+2C_1\log^{+}\frac{\rho}{\rho-r}
 +C_1\left(\log^{+}C_2+2\log 3\right)\\
&\leq C_3\left(\log^{+}T(\rho,f)+\log^{+}\log^{+}T(\rho,f)+\log^{+}\frac{\rho}{\rho-r}+1\right)
\end{split}
\end{equation}
for the pair of positive numbers $r$ and $\rho$ such that $1<r_0<r<\rho<2\rho<R<+\infty$ in \eqref{eq2.8X19}. Here $C_3$ is some large positive constant such that
$C_3\geq \max\left\{2C_1, C_1\left(\log^{+}C_2+2\log 2\right)\right\}.$
\end{remark}
\vskip 2mm
\par
The following result due to Zhang \cite{Zhang1999} that improved Li-Yang \cite[Lemma 7]{LiYang1998}:
\vskip 2mm
\par
\begin{lemma}\rm{}(\cite[Lemma 6]{Zhang1999}) \label{Zhang1999}
Let $f_1$ and $f_2$ be two non-constant meromorphic functions such that
$\overline{N}(r,f_j)+\overline{N}\left(r,\frac{1}{f_j}\right)=S(r)$ for $j\in\{1,2\}.$
Then, either $\overline{N}_0(r,1;f_1,f_2)=S(r),$ or there exist two integers $s$ and $t$ with $|s|+|t|>0,$ such that $f^s_1f^t_2=1.$ Here and in what follows, $\overline{N}_0(r,1;f_1,f_2)$ denotes the reduced counting function of the common $1$-points in $|z|<r,$ and $S(r)$ is any quantity such that $S(r)=o(T(r)),$ as $r\not\in E$ and $r\rightarrow\infty,$ where $T(r)=T(r,f_1)+T(r,f_2),$ and $E\subset (0,+\infty)$ denotes a set of finite linear measure.
\end{lemma}
\vskip 2mm
\par For convenience in stating the following results, we follow Rubel-Taylor \cite[pp.55-58]{RubelTaylor1968} to introduce some notations and some notions as follows: suppose that $Z = \{z_n\}_{n=1}^{+\infty}\subset\Bbb{C}\setminus\{0\}$ is an infinite sequence of non-zero complex numbers such that $|z_n|\rightarrow+\infty,$ as $n\rightarrow+\infty.$ Here we mention that each point in the sequence $Z$ is repeated as many times as its multiplicity. Next we  let
\begin{equation}\nonumber
n(r,Z)=\sum\limits_{|z_n|\leq r}1 \, \, \text{and} \, \,  N\left(r,Z\right)=\int_{0}^{r}\frac{n(t,Z)}{t}dt.
\end{equation}
Throughout this paper, we denote by $\lambda=:\lambda(r)$ a growth function that is a positive, non-decreasing, continuous and unbounded function for $0 < r < +\infty$ (cf.\cite[Definition 1.7]{RubelTaylor1968}). We say that a meromorphic function $f$ is of finite $\lambda$-type to mean that there exist some positive constants $A$ and $B$ such that
  \begin{equation}\nonumber
  T(r, f)\leq A\lambda (Br) \, \, \text{ for}  \, \, r > 0,
\end{equation}
where $T(r,f)$ is the Nevanlinna characteristic of the meromorphic function $f$ for $r>0.$ An entire function $f$ will be of finite $\lambda$-type
if and only if there exist some positive constants $A$ and $B$ such that $|f(z)| \leq \exp (A \lambda (B |z|))$ for each complex number $z\in\Bbb{C}$
(cf.\cite[p.55]{RubelTaylor1968}). we say that the sequence $Z$ has finite $\lambda$-density to mean that there exist some positive constants $A$ and $B$ such that
 \begin{equation}\nonumber
 N(r,Z)\leq A\lambda(Br) \, \, \text{for all} \, \, r > 0.
 \end{equation}
(cf.\cite[Definition 1.8]{RubelTaylor1968}). If $Z$ has finite $\lambda$-density, then there are some positive constants $A$ and $B$ such that
$n(r, Z)\leq A\lambda(Br)$ for all $r > 0$ (cf.\cite[p.58, Proposition 1.9]{RubelTaylor1968}). We say that the sequence $Z$ is $\lambda$-balanced to mean that there exist some positive constants $A$ and $B$ such that
\begin{equation}\nonumber
S(r_1,r_2;k:Z)\leq \frac{A\lambda(Br_1)}{r^k_1}+\frac{A\lambda(Br_2)}{r^k_2}
\end{equation}
for all $r_1>0,$ $r_2 > 0$ and each positive integer $k\in\Bbb{Z}^{+}$ (cf.\cite[Definition 1.10]{RubelTaylor1968}). Here and in what follows
$ S(r_1,r_2;k:Z)=S(r_2;k:Z)-S(r_1;k:Z)$
for all $r_1>0,$ $r_2 > 0$ and each positive integer $k\in\Bbb{Z}^{+},$ where
\begin{equation}\nonumber
S(r;k:Z)=\frac{1}{k} \sum\limits_{0<|z_n|\leq r}\left(\frac{1}{z_n}\right)^k
\end{equation}
for all $r>0$ and each positive integer $k\in\Bbb{Z}^{+}$ (cf.\cite[Definitions 1.5 and 1.6]{RubelTaylor1968}). We recall the following definition:
\begin{definition}\rm{}(\cite[p.60, Definition 1.15]{RubelTaylor1968})\label{Definition2.1}
We say that the sequence $Z$ is $\lambda$-admissible to mean that the sequence $Z$ has finite $\lambda$-density and is $\lambda$-balanced.
\end{definition}
\vskip 2mm
\par Based upon the notations and notions mentioned above, we give the following two results that play an important role in proving the main result of this paper. Here we mention that Lemma \ref{Miles1972} proved in \cite{RubelTaylor1968} originally can be found in \cite{Miles1972}:
\vskip 2mm
\par
\begin{lemma}\rm{}(\cite[ p.372, the main lemma]{Miles1972}) \label{Miles1972}
Suppose $Z = \{z_n\}_{n=1}^{+\infty}\subset\Bbb{C}\setminus\{0\}$ is a sequence of non-zero complex numbers with due regard to multiplicities. If $\lambda(r) = \max\{1,N(r,Z)\}$ for $r>0,$ then there exist some absolute positive constants $A$ and $B,$ and there exists a sequence $\tilde{Z} =\{\tilde{z}_n\}_{n=1}^{+\infty}$ containing the sequence $Z$ with due regard to multiplicities, such that (i) $N(r, \tilde{Z}) \leq 5\lambda (4r)$ for $r > 0$ and (ii) for each positive integer $j$ and any pair of positive numbers $\rho$ and $r$ satisfying $\rho>r>0,$ we have
\begin{equation}\nonumber
S(r,\rho;j:Z)=\left|\frac{1}{j}\sum\limits_{r<|\tilde{z}_n|\leq \rho}\left(\frac{1}{\tilde{z}_n}\right)^j\right|\leq \frac{A\lambda(Br)}{r^j}+\frac{A\lambda(B\rho)}{\rho
^j}.
\end{equation}
\end{lemma}
\vskip 2mm
\par
\begin{lemma}\rm{}(\cite[p.60, Theorem 5.2]{RubelTaylor1968})\label{RubelTaylor1968}
A necessary and sufficient condition that the sequence $Z= \{z_n\}_{n=1}^{+\infty}\subset\Bbb{C}\setminus\{0\}$ be the precise sequence of zeros of an entire function $f$ of finite $\lambda$-type is that the sequence $Z$ be $\lambda$-admissible in the sense of Definition \ref{Definition2.1}, that is to say, the sequence $Z$ have finite $\lambda$-density and be $\lambda$-balanced.
\end{lemma}
\vskip 2mm
\par The following result found in \cite{Conway1978} is called Open Mapping Theorem:
\vskip 2mm
\par
\begin{lemma}\rm{} (\cite[p.99, Open Mapping Theorem 7.5]{Conway1978})\label{Conway1978p.99}
Let $G\subseteq \Bbb{C}$ be a region, and let $\varphi$ be a non-constant analytic function on $G.$ Then for any open set $U$ in $G,$ $\varphi( U)$ is an open set.
\end{lemma}
\vskip 2mm
\par Following Conway \cite[p.39, Definition 2.18 and Proposition 2.19 ]{Conway1978}, we introduce the following notion and one of its propositions: if $G\subseteq\Bbb{C}\setminus \{0\}$ is an open and connected set in $\Bbb{C}$ and $f: G \subseteq\Bbb{C}$ is a continuous function such that $z = \exp {f(z)}$ for all $z$ in $G,$ then $f$ is a branch of the logarithm. If $G \subseteq \Bbb{C}\setminus\{0\}$ is an open and connected set, and $f$ is a branch of $\log z$ on $G,$ then the totality of branches of $\log z$ are the functions $f(z) + 2\pi ki$ with $k\in \Bbb{Z}.$ Next we follow
the proof of  \cite[p.1, Theorem 1.1]{Hayman1964} and the proof of \cite[p.1, Theorem 1.1]{Yang1993} to get the following result that plays an important role in proving the main result of this paper:
\vskip 2mm
\par
\begin{lemma}\rm{} \label{LiWangYi20261}
Let $f$ be a non-constant meromorphic function on the closed disc $|z|\leq R$ with $0<R<+\infty$ such that $f(0)\not\in \{0,\infty\},$
and let $\log f$ be a single-valued analytic branch of the logarithms of $f$ in an open and connected region $G=\Bbb{C}\setminus L,$ where $L$ is a generalized simple directed curve such that its end point is $\infty,$ and its starting point is the point belonging to the set consisting of all the zeros and poles of $f$ in the complex plane such that the starting point has the shortest distance from the origin point of all the points of the set consisting of all the zeros and poles of $f$ in the complex plane, and such that all the zeros and poles of $f$ in the complex plane are on $L.$
Suppose that $\left\{a_{n}\right\}_{n=1}^{M}$ is the multi-set of the no-zero zeros of $\log f$ in the open disc $|z|< R,$ where each point in $\left\{a_{n}\right\}_{n=1}^{M}$ is counted according to its multiplicity of a zero of $\log f.$ Then we have
\begin{equation}\label{eq2.12X19}
m\left(R,\log f\right)=T\left(R,\frac{1}{\log f}\right)+\log |c_{\tau_0}|,
\end{equation}
Here $c_{\tau_0}\neq 0$ is the coefficient of the first term of the Taylor series
\begin{equation}\label{eq2.13X19}
\log f(z)=\sum\limits_{j=\tau_0}^{+\infty} c_jz^j \, \,  \text{with} \, \,  z\in\{z\in\Bbb{C}:|z|<\delta\}
\end{equation}
of $\log f(z)$ at the origin point $z=0,$ where $\delta$ is some small positive constant.
\end{lemma}
\begin{remark}\rm{} \label{remark2.5}
Based upon the desired assertion of Lemma \ref{LiWangYi20261} we introduce the following notations and notions: let $f$ be a non-constant meromorphic function in the complex plane, and let $\left\{a_{n}\right\}_{n=1}^{M}$ be the multi-set of the no-zero zeros of $\log f$ in the open disc $|z|< R$ with $0<R<+\infty,$ where each point in $\left\{a_{n}\right\}_{n=1}^{M}$ is repeated as many times as its multiplicity of a zero of $\log f.$ For $0<r\leq R<+\infty$ we define $m\left(r,\log f\right),$ $m\left(r,\frac{1}{\log f}\right)$ and $N\left(r,\frac{1}{\log f}\right)$ as
\begin{equation}\label{eq2.14X19R1}
m(r,\log f)=\frac{1}{2\pi}\int_0^{2\pi}\log^{+}\left|\log f(re^{i\theta})\right|d\theta,
 \end{equation}
  \begin{equation}\label{eq2.15X19R1}
 m\left(r,\frac{1}{\log f}\right)=\frac{1}{2\pi}\int_0^{2\pi}\log^{+}\left|\frac{1}{\log f(re^{i\theta})}\right|d\theta
 \end{equation}
 and
 \begin{equation}\label{eq2.16X19R1}
N\left(r,\frac{1}{\log f}\right)=\sum\limits_{n=1}^{M}\log \left|\frac{r}{a_n}\right|+n\left(0,\frac{1}{\log f}\right)\log r
 \end{equation}
respectively. Here and in what follows, $n\left(0,\frac{1}{\log f}\right)$ denotes the multiplicity of the origin point of a zero of $\log f.$
We call $m\left(r,\log f\right)$ and $m\left(r,\frac{1}{\log f}\right)$ defined as \eqref{eq2.14X19R1} and \eqref{eq2.15X19R1} the proximity functions of $\log f$ and $\frac{1}{\log f}$ on the circle $|\zeta|=R$ respectively. We call $N\left(r,\frac{1}{\log f}\right)$ defined as \eqref{eq2.16X19R1} the counting function of zeros of $\log f$ in the open disc $|\zeta|<r.$ For $r>0$ we also define
\begin{equation}\label{eq2.17X19R27}
T\left(r,\frac{1}{\log f}\right)=m\left(r,\frac{1}{\log f}\right)+N\left(r,\frac{1}{\log f}\right).
\end{equation}
We call $T\left(r,\frac{1}{\log f}\right)$ the characteristic function of $\frac{1}{\log f}.$
\end{remark}
\vskip 2mm
\par We give the detailed proof of Lemma \ref{LiWangYi20261} as follows:
\vskip 2mm
\par
\begin{proof}
First of all, it is obvious that the zeros and the poles of $f$ on the closed disc $|z|\leq R$ are isolated singularities of the function $\log f$ that are not poles of $\log f.$ Moreover, if we suppose that $\{\tilde{a}_{l}\}_{l=1}^{\tilde{M}_1}$ is the multi-set of zeros of $f$ in the closed disc $|z|\leq R,$ that $\{\tilde{b}_{m}\}_{m=1}^{\tilde{N}_1}$ is the multi-set of poles of $f$ in the closed disc $|z|\leq R,$ and that $\{a_n\}_{n=1}^{M_1}$ is the multi-set of zeros of $\log f$ in the closed disc $|z|\leq R,$ where and in what follows, each point in $\{\tilde{a}_{l}\}_{l=1}^{\tilde{M}_1}$ is counted according to its multiplicity of a zero of $f,$ each point in $\{\tilde{b}_{m}\}_{m=1}^{\tilde{N}_1}$ is counted according to its multiplicity of a pole of $f,$ and each point in $\{a_n\}_{n=1}^{M_1}$ is counted according to its multiplicity of a zero of $\log f,$ then we have
\begin{equation}\label{eq2.17X19R1}
\lim\limits_{z\rightarrow \tilde{a}_{l} }\left|\log f(z)\right|=+\infty \, \, \text{and} \, \, \lim\limits_{z\rightarrow \tilde{b}_{m} }\left|\log f(z)\right|=+\infty,
\end{equation}
\begin{equation}\label{eq2.18X19}
\left|\int_{|z-\tilde{a}_{l}|=\varepsilon}\log\log f(z)dz\right| +\left|\int_{|z-\tilde{b}_{m}|=\varepsilon}\log\log f(z)dz\right|=O\left(\varepsilon\log\log\frac{1}{\varepsilon}\right)
\end{equation}
and
\begin{equation}\label{eq2.19X29}
\left|\int_{|z-a_{n}|=\varepsilon}\log\log f(z)dz\right|=O\left(\varepsilon\log\frac{1}{\varepsilon}\right)
\end{equation}
for the small positive number $\varepsilon$ and any three positive integers $l,$ $m$ and $n$ satisfying $1\leq l\leq \tilde{M}_1,$ $1\leq m\leq \tilde{N}_1$ and $1\leq n\leq M_1$ respectively. Without loss of generality, we next denote by $L_r=:\{z\in L: |z|\leq r\}$ the simple directed curve segment as a part of the directed curve $L$ for $0<r\leq R,$ such that $L_r$ and $L$ have the same starting point in $|z|<r,$ and the end point of $L_r$ is on the circle $|z|=r$ such that all the points except the end point of $L_r$ on $L_r$ belong to the open disc $|z|<r.$ We denote by $L^{-1}_r$ the oppositely directed curve segment of $L_r.$ We denote by
$D_{r, L_r, L^{-1}_r}\subset\Bbb{C}$ the region such that
\begin{equation}\label{eq2.20X25}
D_{r, L_r, L^{-1}_r}=\left\{z\in\Bbb{C}: |z|<r \, \, \text{and} \, \, z\not\in L_r\cup L^{-1}_r\right\}
\end{equation}
with $C_r+L^{-1}_r+L_r$ being the positive boundary of $D_{r, L_r,L^{-1}_r}$ for $0<r\leq R,$ where and in what follows
\begin{equation}\label{eq2.21X25}
C_r=\left\{\zeta\in\Bbb{C}:\zeta=re^{i\theta} \, \text{ with the real variable}  \, \, \theta \, \, \text{from}  \, \, \theta(r) \, \text{ to} \, \, \theta(r)+2\pi\right \}
\end{equation}
denotes a circle with positive orientation such that its radius is equal to $r$ and its center is the origin point for $0<r\leq R.$ Here and in what follows, $\theta(r)\in [-\pi,\pi)$ denotes the principal argument of the endpoint of $L_r$ for $0<r\leq R.$ We consider the following three cases:
\vskip 2mm
\par {\bf Case 1.} Suppose that $f$ has no zeros or poles on the closed disc $|z|\leq R.$ Then $\log f$ is a holomorphic function on the closed
region
 \begin{equation}\label{eq2.22X25}
 \overline{D}_{R, L_R, L^{-1}_R}=\left\{z\in\Bbb{C}: |z|\leq R \, \, \text{and} \, \, z\not\in L_R\cup L^{-1}_R\right\}.
\end{equation}
Next we use \eqref{eq2.20X25}-\eqref{eq2.22X25} and the lines of reasoning of Hayman \cite[p.2, Case (i)]{Hayman1964} for proving Hayman \cite[p.2, (1.3)]{Hayman1964} to get
\begin{equation}\label{eq2.24X25}
\begin{split}
 \log\log f(0)&=\frac{1}{2\pi}\int_{\theta(R)}^{\theta(R)+2\pi}\log\log f(Re^{\phi})d\phi+\frac{1}{2\pi i}\int_{L^{-1}(R)}\frac{\log\log f(\zeta)}{\zeta}d\zeta\\
&\quad+\frac{1}{2\pi i}\int_{L(R)}\frac{\log\log f(\zeta)}{\zeta}d\zeta
\end{split}
 \end{equation}
 with
\begin{equation}\label{eq2.25X25}
\text{Re}\left(\frac{1}{2\pi i}\int_{L^{-1}(R)}\frac{\log\log f(\zeta)}{\zeta}d\zeta\right)=-\text{Re}\left(\frac{1}{2\pi i}\int_{L(R)}\frac{\log\log f(\zeta)}{\zeta}d\zeta\right).
 \end{equation}
We take real parts of both sides of \eqref{eq2.24X25}, and use \eqref{eq2.25X25} and \eqref{eq2.14X19R1}-\eqref{eq2.17X19R27} for $r=R$ to get the desired assertion \eqref{eq2.12X19} with $N\left(R,\frac{1}{\log f}\right)=0$ and $c_{\tau_0}=\log f(0).$
\vskip 2mm
\par {\bf Case 2.} Suppose that $f$ has a finite number of zeros and poles on the circle $|z|=R,$ and that $\log f$ has a finite number of zeros on the circle $|z|=R,$ but $f$ has no zeros and poles in the open disc $|z|<R,$ and $\log f$ has no zeros in the open disc $|z|<R.$ Without loss of generality, we suppose that $\{\tilde{a}_{l}\}_{l=1}^{\tilde{M}_1}$ is the multi-set of zeros of $f$ on the circle $|z|=R,$ that $\{\tilde{b}_{m}\}_{m=1}^{\tilde{N}_1}$ is the multi-set of poles of $f$ on the circle $|z|=R,$ and that $\{a_n\}_{n=1}^{M_1}$ is the multi-set of zeros of $\log f$ on the circle $|z|=R.$ Then we have
\begin{equation}\label{eq2.26X29}
|\tilde{a}_{l}|=|\tilde{b}_{m}|=|a_n|=R
\end{equation}
for any three positive integers $l,$ $m$ and $n$ satisfying $1\leq l\leq \tilde{M}_1,$ $1\leq m\leq \tilde{N}_1$ and $1\leq n\leq M_1$ respectively. Next we denote by $D_{1,R,\varepsilon}$ the region that consists of those points of the open disc $|z|< R$ of the distance greater than the small positive number $\varepsilon$ from every distinct point in the set $\{\tilde{a}_{l}\}_{l=1}^{\tilde{M}_1}\cup \{\tilde{b}_{m}\}_{m=1}^{\tilde{N}_1}\cup\{a_n\}_{n=1}^{M_1}.$ We denote by $L_{1,R,\varepsilon}$ the boundary of the region $D_{1,R,\varepsilon}.$ Then $\log f$ is a holomorphic function on the closed region $\overline{D}_{1,R,\varepsilon}=D_{1,R,\varepsilon}\cup L_{1,R,\varepsilon}.$ Next we let
\begin{equation}\label{eq2.21X19}
g_1(\zeta)=
\begin{cases}
\left(\frac{R}{\zeta}\right)^{\tau_0}\log f(\zeta), &\text{if} \quad \zeta\in \overline{D}_{1,R,\varepsilon}\setminus\{0\},\\
c_{\tau_0}R^{\tau_0},&\text{if}\quad  \zeta=0.\\
\end{cases}
\end{equation}
Noting that $c_{\tau_0}\neq 0$ is the coefficient of the first term of the Taylor series \eqref{eq2.13X19}, we have from \eqref{eq2.13X19} that the function $g_1(\zeta)$ defined as \eqref{eq2.21X19} is continuous at the origin point $\zeta=0.$ Combining this with Ahlfors \cite[p.122, Morera's theorem]{Ahlfors1979}, we deduce that the function $g_1(\zeta)$ defined as \eqref{eq2.21X19} is analytic at the origin point $\zeta=0.$ Therefore, the function $(\log g_1(\zeta))/\zeta$ is holomorphic on the closed region $\overline{D}_{1,R,\varepsilon}=D_{1,R,\varepsilon}\cup L_{1,R,\varepsilon},$ except for the point $\zeta=0,$ which is possibly a pole with residue $\log g_1(0).$ Therefore, we have
\begin{equation}\label{eq2.22X19}
\log g_1(0)=\frac{1}{2\pi i} \int_{L_{1,R,\varepsilon}}\frac{\log g_1(\zeta)}{\zeta}d\zeta.
\end{equation}
Since the length of every small circular arcs included in $L_{1,R,\varepsilon}$ of radius $\varepsilon$ and centred on the points from the set $\{\tilde{a}_{l}\}_{l=1}^{\tilde{M}_1}\cup \{\tilde{b}_{m}\}_{m=1}^{\tilde{N}_1}\cup\{a_n\}_{n=1}^{M_1}$ is not greater than $2\pi\varepsilon,$ we have that the total length of these small circular arcs is $O(\varepsilon),$ the absolute value of the integrand along those small circular arcs centred on the points from the set $\{\tilde{a}_{l}\}_{l=1}^{\tilde{M}_1}\cup \{\tilde{b}_{m}\}_{m=1}^{\tilde{N}_1}$ is $O\left(\log\log\frac{1}{\varepsilon}\right)$ uniformly, and the absolute value of the integrand along those small circular arcs centred on the points from the set $\{a_n\}_{n=1}^{M_1}$ is $O\left(\log\frac{1}{\varepsilon}\right)$ uniformly. This together with \eqref{eq2.17X19R1}-\eqref{eq2.19X29} and \eqref{eq2.26X29} implies that the sum of the absolute value of the integrals along the small circular arcs centred on the points from the set $\{\tilde{a}_{l}\}_{l=1}^{\tilde{M}_1}\cup \{\tilde{b}_{m}\}_{m=1}^{\tilde{N}_1}\cup\{a_n\}_{n=1}^{M_1}$ is not greater than $O\left(\varepsilon\log\log\frac{1}{\varepsilon}\right)+O\left(\varepsilon\log\frac{1}{\varepsilon}\right).$ Combining this with making $\varepsilon\rightarrow 0^{+}$ on both sides of \eqref{eq2.22X19}, we have
\begin{equation}\label{eq2.23X19}
\log g_1(0)=\frac{1}{2\pi i} \int_{L_{1,R}}\frac{\log g_1(\zeta)}{\zeta}d\zeta,
\end{equation}
where and in what follows
\begin{equation}\label{eq2.24X19}
L_{1,R}=\left\{\zeta\in\Bbb{C}:\zeta=Re^{i\theta} \, \text{ with the real variable}  \, \, \theta \, \, \text{from}  \, \, 0  \, \text{ to} \, 2\pi\right \}
\end{equation}
denotes a circle with positive orientation such that its radius is equal to $R$ and its center is the origin point. Noting that $\left|\left(\frac{R}{\zeta}\right)^{\tau_0}\right|=1$ for each point $\zeta$ such that $|\zeta|=R,$ we substitute \eqref{eq2.21X19} into \eqref{eq2.23X19} and then take the real parts on both sides of the obtained equality to have
\begin{equation}\label{eq2.25X20}
\log |c_{\tau_0}|+\tau_0 \log R=\frac{1}{2\pi}\int_0^{2\pi}\log\left|\log f(Re^{i\theta})\right|d\theta.
\end{equation}
On the other hand, from \eqref{eq2.13X19} and the supposition at the beginning of Case 2 we have $\tau_0=0$ and
$N\left(R,\frac{1}{\log f}\right)=0.$ Combining this with \eqref{eq2.25X20} and
\begin{equation}\label{eq2.26X20}
\begin{split}
&\quad\frac{1}{2\pi}\int_0^{2\pi}\log\left|\log f(Re^{i\theta})\right|d\theta=\frac{1}{2\pi}\int_0^{2\pi}\log^{+}\left|\log f(Re^{i\theta})\right|d\theta\\
&\quad -\frac{1}{2\pi}\int_0^{2\pi}\log^{+}\left|\frac{1}{\log f(Re^{i\theta})}\right|d\theta
=m\left(R,\log f\right)-m\left(R,\frac{1}{\log f}\right),
\end{split}
\end{equation}
we have
\begin{equation}\nonumber
m\left(R,\log f\right)=m\left(R,\frac{1}{\log f}\right)+\log |c_{\tau_0}|=T\left(R,\frac{1}{\log f}\right)+\log |c_{\tau_0}|
\end{equation}
 This is a special case of the desired assertion \eqref{eq2.12X19}.
\vskip 2mm
\par {\bf Case 3.} Suppose that $f$ has a finite number of zeros and poles in the open disc $|z|<R,$ and that $\log f$ has a finite number of zeros in the open disc $|z|<R.$ Without loss of generality, we suppose that $\{\tilde{a}_{l}\}_{l=1}^{\tilde{M}_1}$ is the multi-set of zeros of $f$ in the open disc $|z|< R,$ that $\{\tilde{b}_{m}\}_{m=1}^{\tilde{N}_1}$ is the multi-set of poles of $f$ in the open disc $|z|< R,$ and that $\{a_n\}_{n=1}^{M_1}$ is the multi-set of zeros of $\log f$ in the open disc $|z|< R.$ Without loss of generality, we also suppose that the set $\{\tilde{a}_{l}\}_{l=1}^{\tilde{M}_1}$ has $\tilde{M}_1$ distinct elements, and the set $\{\tilde{b}_{m}\}_{m=1}^{\tilde{N}_1}$ has $\tilde{N}_1$ distinct elements. Next we set
\begin{equation}\label{eq2.25X19}
\gamma^{-1}_{1,l,\varepsilon}:\zeta-\tilde{a}_{l}=\varepsilon e^{i\theta} \, \, \text{with the real variable} \, \, \theta  \, \, \text{from} \, \, 2\pi  \, \, \text{to} \, \,  0
\end{equation}
for the positive integer $l$ satisfying $1\leq l\leq \tilde{M}_1,$ and set
\begin{equation}\label{eq2.26X19}
\gamma^{-1}_{2,m,\varepsilon}:\zeta-\tilde{b}_{m}=\varepsilon e^{i\theta} \, \, \text{with the real variable}  \, \, \theta  \, \, \text{from} \, \, 2\pi  \, \, \text{to} \, \,  0
\end{equation}
for the positive integer $m$ satisfying $1\leq m\leq \tilde{N}_1.$ Here the positive number $\varepsilon$ is small enough such that the $\tilde{m}_1+\tilde{N}_1$ distinct circles
\begin{equation}\nonumber
\gamma^{-1}_{1,1,\varepsilon}, \gamma^{-1}_{1,2,\varepsilon}, \ldots, \gamma^{-1}_{1,\tilde{M}_1,\varepsilon}, \gamma^{-1}_{2,1,\varepsilon}, \gamma^{-1}_{2,2,\varepsilon}, \ldots, \gamma^{-1}_{2,\tilde{N}_1,\varepsilon}
\end{equation}
are in the open disc $|z|<R$ and separated with each other, and that each point in the set  $\{a_n\}_{n=1}^{M_1}$ is outside of all these $\tilde{m}_1+\tilde{N}_1$ distinct circles. Then, the function
\begin{equation}\label{eq2.27X19}
g_2(\zeta)=: \frac{g_1(\zeta)}{\prod\limits_{n=1}^{M_1}\frac{R(\zeta-a_n)}{R^2-\overline{a}_n\zeta}}
\end{equation}
is holomorphic in the region $D_{2,R,\varepsilon},$ where $g_1(\zeta)$ is the function defined as \eqref{eq2.21X19} and $D_{2,R,\varepsilon}$ denotes the region that consists of those points of the open disc $|z|< R$ of the distance greater than the small positive number $\varepsilon$ from every distinct point in the set $\{\tilde{a}_{l}\}_{l=1}^{\tilde{M}_1}\cup \{\tilde{b}_{m}\}_{m=1}^{\tilde{N}_1}.$ In addition, from \eqref{eq2.25X19}, \eqref{eq2.26X19} and the assertion of Case 2 we suppose, without loss of generality, that
\begin{equation}\label{eq2.28X19}
L_{2,R,\varepsilon}=L_{1,R}+L^{-1}_R+L_R+\gamma^{-1}_{1,1,\varepsilon}+\gamma^{-1}_{1,2,\varepsilon}+\cdots+\gamma^{-1}_{1,\tilde{M}_1,\varepsilon}
+\gamma^{-1}_{2,1,\varepsilon}+\gamma^{-1}_{2,2,\varepsilon}+ \cdots+\gamma^{-1}_{2,\tilde{N}_1,\varepsilon}
\end{equation}
is the boundary of the region $D_{2,R,\varepsilon}$ with positive orientation. Therefore, from Cauchy's residue theorem and the similar reasoning in Case 2 we have
\begin{equation}\label{eq2.29X19}
\log g_2(0)=\frac{1}{2\pi i} \int_{L_{2,R,\varepsilon}}\frac{\log g_2(\zeta)}{\zeta}d\zeta,
\end{equation}
On the other hand, from the assumption $f(0)\not\in\{0,\infty\}$ we have $0\not\in\{\tilde{a}_{l}\}_{l=1}^{\tilde{M}_1}\cup \{\tilde{b}_{m}\}_{m=1}^{\tilde{N}_1}.$ Combining this with \eqref{eq2.21X19}, \eqref{eq2.25X19}, \eqref{eq2.26X19} and \eqref{eq2.27X19}, we can verify that
\begin{equation}\label{eq2.30X19R1}
\begin{split}
&\quad\left|\int_{\gamma^{-1}_{1,l,\varepsilon}}  \frac{\log g_2(\zeta)}{\zeta}d\zeta\right|+\left|\int_{\gamma^{-1}_{2,m,\varepsilon}}  \frac{\log g_2(\zeta)}{\zeta}d\zeta\right|
=\left|\int_{\gamma^{-1}_{1,l,\varepsilon}}\frac{1}{\zeta}\log \left(\frac{\left(\frac{R}{\zeta}\right)^{\tau_0}\log f(\zeta)}{\prod\limits_{n=1}^{M_1}\frac{R(\zeta-a_n)}{R^2-\overline{a}_n\zeta}}\right)\right|\\
&\quad +\left|\int_{\gamma^{-1}_{2,m,\varepsilon}}\frac{1}{\zeta}\log \left(\frac{\left(\frac{R}{\zeta}\right)^{\tau_0}\log f(\zeta)}{\prod\limits_{n=1}^{M_1}\frac{R(\zeta-a_n)}{R^2-\overline{a}_n\zeta}}\right)\right|=O\left(\varepsilon\log\log\frac{1}{\varepsilon}\right)
\end{split}
\end{equation}
for the small positive number $\varepsilon$ and any two positive integers $l$ and $m$ satisfying $1\leq l\leq \tilde{M}_1$ and $1\leq m\leq \tilde{N}_1$ respectively. Whether one of $f(\zeta),$ $1/f(\zeta)$ and $\log f(\zeta)$ has zeros or not on the circle $|\zeta|=R$ for the positive number $R,$ we apply   \eqref{eq2.25X19}-\eqref{eq2.28X19}, \eqref{eq2.30X19R1} and the line of the reasoning of \eqref{eq2.23X19} in Case 2 to make $\varepsilon\rightarrow 0^{+}$ on both sides of \eqref{eq2.29X19}, and then we have
\begin{equation}\label{eq2.31X19R1}
\log g_2(0)=\frac{1}{2\pi i} \int_{L_{1,R}}\frac{\log g_2(\zeta)}{\zeta}d\zeta+\int_{L^{-1}_R}\frac{\log g_2(\zeta)}{\zeta}d\zeta+\int_{L_R}\frac{\log g_2(\zeta)}{\zeta}d\zeta,
\end{equation}
From \eqref{eq2.24X19} we deduce that \eqref{eq2.31X19R1} can be rewritten into
\begin{equation}\label{eq2.32X19R1}
\log g_2(0)=\frac{1}{2\pi} \int_{0}^{2\pi}\log g_2(Re^{i\theta})d\theta+\int_{L^{-1}_R}\frac{\log g_2(\zeta)}{\zeta}d\zeta+\int_{L_R}\frac{\log g_2(\zeta)}{\zeta}d\zeta
\end{equation}
with
\begin{equation}\label{eq2.41X29}
\text{Re}\left(\frac{1}{2\pi i}\int_{L^{-1}(R)}\frac{\log g_2(\zeta)}{\zeta}d\zeta\right)=-\text{Re}\left(\frac{1}{2\pi i}\int_{L(R)}\frac{\log g_2(\zeta)}{\zeta}d\zeta\right).
 \end{equation}
On the other hand, from \eqref{eq2.21X19} and \eqref{eq2.27X19} we have
\begin{equation}\label{eq2.33X19R1}
g_2(0)=\frac{g_1(0)}{\frac{(-1)^{M_1}}{R^{M_1}}\prod\limits_{n=1}^{M_1}a_n}=\frac{ c_{\tau_0}R^{\tau_0}}{\frac{(-1)^{M_1}}{R^{M_1}}\prod\limits_{n=1}^{M_1}a_n}
\end{equation}
and
\begin{equation}\label{eq2.34X19R1}
g_2(Re^{i\theta})=\frac{g_1(Re^{i\theta})}{\prod\limits_{n=1}^{M_1}\frac{R(Re^{i\theta}-a_n)}{R^2-\overline{a}_nRe^{i\theta}}}=\frac{ \left(\frac{R}{Re^{i\theta}}\right)^{\tau_0}\log f(Re^{i\theta})}{\prod\limits_{n=1}^{M_1}\frac{R(Re^{i\theta}-a_n)}{R^2-\overline{a}_nRe^{i\theta}}}
\end{equation}
for $\zeta=Re^{i\theta}$ with $\theta\in [0,2\pi).$ Noting that
\begin{equation}\nonumber
\left|\left(\frac{R}{\zeta}\right)^{\tau_0}\right|=\left| \prod\limits_{n=1}^{M_1}\frac{R(\zeta-a_n)}{R^2-\overline{a}_n\zeta} \right|=1 \, \, \text{for each point} \, \, \zeta \, \text{such that} \, \, |\zeta|=R,
\end{equation}
Noting that \eqref{eq2.41X29} holds, we substitute \eqref{eq2.33X19R1} and \eqref{eq2.34X19R1} into \eqref{eq2.32X19R1}, and then take the real parts on both sides of the obtained equality to have
\begin{equation}\label{eq2.35X19R1}
\frac{1}{2\pi}\int_0^{2\pi}\log\left|\log f(Re^{i\theta})\right|d\theta=\sum\limits_{n=1}^{M}\log \left|\frac{R}{a_n}\right|+\log |c_{\tau_0}|+\tau_0 \log R.
\end{equation}
and
\begin{equation}\label{eq2.37X19R1}
N\left(R,\frac{1}{\log f}\right)=\sum\limits_{n=1}^{M}\log \left|\frac{R}{a_n}\right|+n\left(0,\frac{1}{\log f}\right)\log R=\sum\limits_{n=1}^{M}\log \left|\frac{R}{a_n}\right|+\tau_0 \log R,
 \end{equation}
we substitute \eqref{eq2.26X20} and \eqref{eq2.37X19R1} into \eqref{eq2.35X19R1}, and then we have
\begin{equation}\nonumber
\begin{split}
m\left(R,\log f\right)&=m\left(R,\frac{1}{\log f}\right)+N\left(R,\frac{1}{\log f}\right)+\log |c_{\tau_0}|=T\left(R,\frac{1}{\log f}\right)+\log |c_{\tau_0}|.
\end{split}
\end{equation}
This is the desired assertion \eqref{eq2.12X19} of Lemma \ref{LiWangYi20261}.
\end{proof}
\par
\vskip 2mm
The following result is called Chi-Tai Chuang$'$ inequality proved by Chuang\cite{Chuang1951} originally was also proved in Yang \cite[pp.95-96]{Yang1993}:
\par
\vskip 2mm
\begin{lemma}\rm{}(\cite{Chuang1951})\label{Chuang1951}
Suppose that $f$ is a non-constant meromorphic function in the complex plane. If $f(0)\neq \infty,$ then
$T(r,f)<C_{\tau}T(\tau r,f')+\log^{+}(\tau r)+4+\log^{+}|f(0)|$
for $\tau>1$ and $r>0.$ Here $C_{\tau}$ is a positive constant depending only upon the positive constant $\tau.$
\end{lemma}
\par
\vskip 2mm
The following result is due to Yang \cite{Yang1993}:
\par
\vskip 2mm
\begin{lemma}\rm{}(\cite[p.47, Lemma 2.3]{Yang1993})\label{Yang1993}
For a real number $x>0$ and a positive constant $A$ satisfying $A\geq e,$ we have
$\log x+A\log^{+}\log^{+}\frac{1}{x}\leq \log^{+}x+A\left(\log A-1\right).$
\end{lemma}
\vskip 2mm
\par The following result is due to Ku\cite{Ku1978}:
\par
\vskip 2mm
\begin{lemma}\rm{}(\cite[Lemma 2]{Ku1978})\label{Ku1978}
Let $T(r): (0,R)\longrightarrow [0,+\infty)$ be a continuous, non-decreasing and non-negative function, where $0<R\leq +\infty,$ and let $A,$ $B$ and $C$ be three positive constants. If
\begin{equation}\nonumber
T(r)<A\log^{+}T(\rho)+B\left(\log^{+}\frac{1}{\rho-r}+\log^{+}\rho+\log^{+}\frac{1}{r}\right)+C
\end{equation}
for $0<r<\rho<R<+\infty.$ Then, we have
\begin{equation}\nonumber
T(r)<4(A+B)\left(\log\frac{\rho}{\rho-r}+\log^{+}R+\log^{+}\frac{1}{R}\right)+20(A+B+1)^2+2C
\end{equation}
 for $0<r<\rho<R<+\infty.$
\end{lemma}
\vskip 2mm
\par
	\section{Proof of Theorem \ref{Theorem1.2}}
\vskip 2mm
\par
First of all, we introduce some notations as follows: let $f_1$ and $f_2$ be two non-constant meromorphic functions in the complex plane, and let $a\in \mathbb{C}$ be a finite complex value. Next we denote by $\overline{N}_{(1,l)}(r,a;f_1,f_2)$ the reduced counting function of those common zeros of $f_1-a$ and $f_2-a$ in $\left|z\right|<r,$ where each point in $\overline{N}_{(1,l)}(r,a;f_1,f_2)$ is a simple zero of $f_1-a$ and a zero of $f_2-a$ of multiplicity $l$ for the given positive integer $l.$ We denote by $\overline{N}_{(m,1)}(r,a;f_1,f_2)$ the reduced counting function of those common zeros of $f_1-a$ and $f_2-a$ in $\left|z\right|<r,$ where each point in $\overline{N}_{(m,1)}(r,a;f_1,f_2)$ is a simple zero of $f_2-a$ and a zero of $f_1-a$ of multiplicity $m$ for the given positive integer $m.$
\vskip 2mm
\par Secondly, we want to point out that the lines from the next formal beginning of the proof of Theorem \ref{Theorem1.2} to the end of Case 1 of the next formal beginning of the proof of Theorem \ref{Theorem1.2} are covered essentially by the proof of Steinmetz \cite[Theorem 6.2]{Steinmetz2012} without any assumption on the order of growth. For the convenience of the readers, we give the detailed proof for this part without any assumption on the order of growth.
\vskip 2mm
\par Thirdly, we prove Theorem \ref{Theorem1.2} formally: first of all, we use Lemma \ref{Adams1971} and the assumption that $f$ and $g$ are two distinct non-constant meromorphic functions that share three distinct finite values $a_1,a_2,a_3$ IM and share $\infty$ CM to deduce that $f$ and $g$ are two distinct transcendental meromorphic functions. Without loss of generality, we suppose that
 \begin{equation}\label{eq3.1X19}
 \left\{ f(0)-a, g(0)-a,f'(0), g'(0)\right\}\cap \{0,\infty\}=\emptyset
 \end{equation}
 and
  \begin{equation}\label{eq3.2X19}
\left\{\frac {lf'(0)\prod\limits_{j=1}^3(g(0)-a_j)} {g'(0)\prod\limits_{j=1}^3(f(0)-a_j)}, \frac {f'(0)\prod\limits_{j=1}^3(g(0)-a_j)} {mg'(0)\prod\limits_{j=1}^3(f(0)-a_j)}\right\}\cap \{0,1,\infty\}=\emptyset
 \end{equation}
 for each complex number $a\in\{0,a_1,a_2,a_3\}$ and any pair of positive integers  $l$ and $m$ such that $1\leq l\leq 6$ and $2\leq m\leq 6.$ Here and in what follows, $\emptyset$ denotes the empty set. Otherwise, we replace $f(z)$ and $g(z)$ with $f(z+\eta)$ and $g(z+\eta)$ respectively, such that
 \begin{equation}\nonumber
 \left\{ f(\eta)-a, g(\eta)-a,f'(\eta), g'(\eta)\right\}\cap \{0,\infty\}=\emptyset
 \end{equation}
 and
  \begin{equation}\nonumber
\left\{\frac {lf'(\eta)\prod\limits_{j=1}^3(g(\eta)-a_j)} {g'(\eta)\prod\limits_{j=1}^3(f(\eta)-a_j)}, \frac {f'(\eta)\prod\limits_{j=1}^3(g(\eta)-a_j)} {mg'(\eta)\prod\limits_{j=1}^3(f(\eta)-a_j)}\right\}\cap \{0,1,\infty\}=\emptyset
 \end{equation}
for each complex number $a\in\{0,a_1,a_2,a_3\}$ and any pair of positive integers  $l$ and $m$ satisfying $1\leq l\leq 6$ and $2\leq m\leq 6,$
where $\eta\neq 0$ is some finite fixed complex number. Based upon this supposition and the assumptions of Theorem \ref{Theorem1.2}, we let $\Phi_{f}$ and $\Phi_{g}$ be defined as \eqref{eq1.1X1}. Then it follows from  Theorem \ref{TheoremG} that
$\Phi_{f}$ and $\Phi_{g}$ are entire functions. Combining this with \eqref{eq1.1X1} and Lemma \ref{Gundersen1992}, we have that the function
\begin{equation}\label{eq3.3X19}
\Phi_{f}\Phi_{g}=\Psi=\frac{f'g'{(f-g)}^2}{(f-a_1)(f-a_2)(f-a_3)(g-a_1)(g-a_2)(g-a_3)}
\end{equation}
 is also an entire function named Mues$'$s function mentioned before Lemma \ref{Gundersen1992}. This together with \eqref{eq2.5X19}, \eqref{eq3.3X19} and \eqref{eq2.11X19} of Remark \ref{remark2.4} gives
\begin{equation}\label{eq3.4X19}
\begin{split}
&\quad N\left(r,\frac{1}{\Phi_{f}}\right)+N\left(r,\frac{1}{\Phi_{g}}\right)=N\left(r,\frac{1}{\Psi}\right)\leq T(r,\Psi)+O(1)=m(r,\Psi)+O(1)\\
&=m\left(r,\sum\limits_{s=1}^3\sum\limits_{t=1}^3\frac{c_{st}f'}{f-a_s}\cdot\frac{g'}{g-a_t}\right)+O(1)\\
&\leq C_3\left(\log^{+}T(\rho,f)+\log^{+}\log^{+}T(\rho,f)+\log^{+}\frac{\rho}{\rho-r}+1\right)
\end{split}
\end{equation}
for any pair of positive numbers $r$ and $\rho$ such that $1<r_0<r<\rho<R$ in \eqref{eq2.7X19}. On the other hand, from Lemma \ref{LiYi2007} we have
\begin{equation}\label{eq3.5X19}
\frac {1}{7} T(r,f)\leq\sum\limits_{j=1}^{3}\sum\limits_{l=1}^{6}\overline{N}_{(1,l)}\left(r,a_j;f,g\right)
+\sum\limits_{j=1}^{3}\sum\limits_{m=2}^{6}\overline{N}_{(m,1)}\left(r,a_j;f,g\right)+S(r,f)+S(r,g).
\end{equation}
From Remarks \ref{remark2.3} and \ref{remark2.4} we see that the error term $S(r,f)+S(r,g)$ in \eqref{eq3.5X19} can be estimated as \eqref{eq2.11X19}. Therefore, the error term $S(r,f)+S(r,g)$ in \eqref{eq3.5X19} can be estimated as
\begin{equation}\label{eq3.6X19}
S(r,f)+S(r,g)\leq C_3\left(\log^{+}T(\rho,f)+\log^{+}\log^{+}T(\rho,f)+\log^{+}\frac{\rho}{\rho-r}+1\right)
\end{equation}
for any pair of positive numbers $r$ and $\rho$ such that $1<r_0<r<\rho<R<+\infty.$ From the expressions of the entire functions $\Phi_{f}$ and $\Phi_{g}$ defined as \eqref{eq1.1X1}, we have
\begin{equation}\label{eq3.7X19}
\frac{\Phi_{f}}{\Phi_{g}}=\frac {f'(g-a_1)(g-a_2)(g-a_3)} {g'(f-a_1)(f-a_2)(f-a_3)}.
\end{equation}
\vskip 2mm
\par
From \eqref{eq3.7X19} we consider the following two cases:
\vskip 2mm
\par {\bf Case 1.} Suppose that $\Phi_{f}/\Phi_{g}$ is reduced to a non-zero constant. Then, it follows from \eqref{eq3.7X19} that
\begin{equation}\label{eq3.8X19}
\frac {f'} {(f-a_1)(f-a_2)(f-a_3)}=\frac {cg'} {(g-a_1)(g-a_2)(g-a_3)}.
\end{equation}
From \eqref{eq3.5X19} and \eqref{eq3.8X19} we consider the following three subcases:
\vskip 2mm
\par {\bf Subcase 1.1.} Suppose that
\begin{equation}\label{eq3.9X19}
\overline{N}_{(1,1)}\left (r,a_1;f,g \right )+\overline{N}_{(1,1)}\left (r,a_2;f,g \right ) +\overline{N}_{(1,1)}\left (r,a_3;f,g \right ) \neq S(r,f).
\end{equation}
Here and in what follows, $\overline{N}_{(1,1)}\left (r,a;f,g \right )$ denotes the reduced counting function of the common simple zeros of $f-a$ and $g-a$ in $|z|<r$ for the value $a\in\{a_1,a_2,a_3\}.$ From \eqref{eq3.9X19} we suppose, without loss of generality, that
\begin{equation}\label{eq3.10X19}
\overline{N}_{(1,1)}\left (r,a_1;f,g \right) \neq S(r,f).
\end{equation}
From \eqref{eq3.10X19} we see that there exists
 a common simple zero $z_{11}\in\Bbb{C} $ of $f-a_1$ and $g-a_1.$ Combining this with \eqref{eq3.8X19} and the assumption that $f$ and $g$ share $a_1,$ $a_2,$ $a_3$ IM and $\infty$ CM, we deduce $c=1,$ and so \eqref{eq3.8X19} can be rewritten into
\begin{equation}\label{eq3.11X19}
\frac {f'} {(f-a_1)(f-a_2)(f-a_3)}=\frac {g'} {(g-a_1)(g-a_2)(g-a_3)}.
\end{equation}
 From \eqref{eq3.11X19} and the assumption that $f$ and $g$ share $a_1,$ $a_1,$ $a_3$ IM and $\infty$ CM, we deduce that $f$ and $g$ share $a_1,$ $a_2,$ $a_3,$ $\infty$  CM. This reveals the assertion of Theorem \ref{Theorem1.2}.
\vskip 2mm
\par {\bf Subcase 1.2.} Suppose that there exists some positive integer $l$ satisfying $2\leq l\leq 6,$ such that
\begin{equation}\label{eq3.12X19}
\overline{N}_{(1,l)}\left (r,a_1;f,g \right )+\overline{N}_{(1,l)}\left (r,a_2;f,g \right ) +\overline{N}_{(1,l)}\left (r,a_3;f,g \right )\neq S(r,f),
\end{equation}
 where and in what follows, $\overline{N}_{(1,l)}\left (r,a;f,g \right )$ denotes the reduced counting function of those common zeros of $f-a$ and $g-a$ in $|z|<r,$ where each such common zero of $f-a$ and $g-a$ is a simple zero of $f-a$ and a zero of $g-a$ of multiplicity $l$ for $a\in\{a_1,a_2,a_3\}.$ From \eqref{eq3.12X19} we suppose, without loss of generality, that
 \begin{equation}\label{eq3.13X19}
\overline{N}_{(1,l)}\left (r,a_1;f,g \right )\neq S(r,f).
\end{equation}
 From \eqref{eq3.13X19} we see that there exists a point $z_{11l}\in \Bbb{C}$  that is a simple zero of $f-a_1,$ and a zero of $g-a_1$ of multiplicity $l.$ Combining this with \eqref{eq3.8X19} and the assumption that $f$ and $g$ share $a_1,$ $a_1,$ $a_3$ IM and $\infty$ CM, we deduce $c=1/l,$ and so \eqref{eq3.8X19} can be rewritten into
\begin{equation}\label{eq3.14X19}
\frac {lf'} {(f-a_1)(f-a_2)(f-a_3)}=\frac {g'} {(g-a_1)(g-a_2)(g-a_3)}.
\end{equation}
 Next we set
   \begin{equation}\label{eq3.15X19}
\varphi_1=\frac {lf'} {(f-a_1)(f-a_2)}-\frac {g'} {(g-a_1)(g-a_2)},
\end{equation}
 \begin{equation}\label{eq3.16X19}
 \varphi_2=\frac {lf'} {(f-a_1)(f-a_3)}-\frac {g'} {(g-a_1)(g-a_3)},
\end{equation}
  \begin{equation}\label{eq3.17X19}
f_1=\frac{(g-a_1)(f-a_2)^l}{(f-a_1)^l(g-a_2)} \, \, \text{and} \, \, f_2=\frac{(g-a_1)(f-a_3)^l}{(f-a_1)^l(g-a_3)}.
\end{equation}
From \eqref{eq3.14X19} and the assumption that $f$ and $g$ share $a_1,$ $a_1,$ $a_3$ IM and $\infty$ CM, we deduce that if $z_{1l}\in\Bbb{C}$ is a zero of $f-a$ of multiplicity $\nu_{1l}$ for $a\in\{a_1,a_2,a_3\},$ then  $z_{1l}\in\Bbb{C}$ is a zero of $g-a$ of multiplicity $l\nu_{1l}\geq l\geq 2.$ This reveals that
\begin{equation}\label{eq3.18X19}
N\left(r,f_1\right)+N\left(r,f_2\right)+N\left(r,\frac{1}{f_1}\right)+N\left(r,\frac{1}{f_2}\right)=0,
\end{equation}
and that the function $\varphi_1$ in \eqref{eq3.15X19} and the function  $\varphi_2$ in \eqref{eq3.16X19} are entire functions.
 Combining this with \eqref{eq3.15X19}, \eqref{eq3.16X19}, Lemma \ref{Nevanlinna1926}(i) and  Yang \cite[p.23,Theorem 1.6]{Yang1993}, we deduce
  \begin{equation}\nonumber 
    T(r,\varphi_1)+ T(r,\varphi_2)=   m(r,\varphi_1)+ m(r,\varphi_2)=S(r,f).
  \end{equation}
 Suppose that $\varphi_1=0.$ Then, it follows from \eqref{eq3.15X19} that
 \begin{equation}\label{eq3.19X19}
 \frac {lf'} {(f-a_1)(f-a_2)}=\frac {g'} {(g-a_1)(g-a_2)}
  \end{equation}
  From \eqref{eq3.14X19} and \eqref{eq3.19X19} we deduce $f=g.$ This contradicts the assumption that $f$ and $g$ are two distinct and non-constant meromorphic functions.  Similarly, we get a contradiction from \eqref{eq3.14X19} and \eqref{eq3.16X19} provided that $\varphi_2=0.$ Therefore, $\varphi_1$ and $\varphi_2$ are non-vanishing entire functions. On the other hand, from the first fundamental theorem, the second fundamental theorem and the known result that the multiplicity of each zero of $g-a$ is greater or equal to $l\geq 2$ for each value $a\in\{a_1,a_2,a_3\},$ we have
\begin{equation}\nonumber
\begin{split}
2T(r,g)&=\overline{N}\left(r,g\right)+\sum\limits_{j=1}^{3}\overline{N}\left(r,\frac{1}{g-a_j}\right)+S(r,g)\leq \overline{N}\left(r,g\right)+N\left(r,\frac{1}{g'}\right)+S(r,g)\\
&=\overline{N}\left(r,g\right)+m(r,g')+N(r,g')-m\left(r,\frac{1}{g'}\right)+S(r,g)\\
&\leq 2\overline{N}\left(r,g\right)+m(r,g)+N(r,g)-m\left(r,\frac{1}{g'}\right)+m\left(r,\frac{g'}{g}\right)+S(r,g)\\
&\leq  2\overline{N}\left(r,g\right)+T(r,g)+S(r,g),
\end{split}
\end{equation}
and so we have
\begin{equation}\label{eq3.20X19}
T(r,g)\leq  2\overline{N}\left(r,g\right)+S(r,g).
\end{equation}
Noting that $f$ and $g$ share $a_1,$ $a_2,$ $a_3$ IM and $\infty$ CM, we have from \eqref{eq3.17X19} and \eqref{eq3.20X19} that each common pole of $f$ and $g$ must be a common zero of $f_1-1$ and $f_2-1.$ Therefore, we have
\begin{equation}\label{eq3.21X19}
\frac{1}{2}T(r,g)\leq \overline{N}_0(r,\infty;f,g)+S(r,g)\leq \overline{N}_0(r,1;f_1,f_2)+S(r,g),
\end{equation}
where and in what follows, $\overline{N}_0(r,\infty;f,g)$ denotes the reduced counting function of the common poles of $f$ and $g$ in $|z|<r,$ and $\overline{N}_0(r,1;f_1,f_2)$ denotes the reduced counting function of the common zeros of $f_1-1$ and $f_2-1$ in $|z|<r.$
From \eqref{eq3.17X19}, Lemma \ref{Nevanlinna1926}(i) and Valiron-Mohon'ko lemma (cf.\cite{Mohon'ko1971}) we have
\begin{equation}\label{eq3.22X19}
\begin{split}
&\quad T(r,f_1)+T(r,f_2)=T\left(r,\frac{(g-a_1)(f-a_2)^l}{(f-a_1)^l(g-a_2)}\right)+T\left(r,\frac{(g-a_1)(f-a_3)^l}{(f-a_1)^l(g-a_3)}\right)\\
                 &\leq T\left(r,\frac{g-a_1}{g-a_2}\right)+T\left(r,\left(\frac{f-a_2}{f-a_1}\right)^l\right)+T\left(r,\frac{g-a_1}{g-a_3}\right)
                 +T\left(r,\left(\frac{f-a_3}{f-a_1}\right)^l\right)\\
                 &=T(r,g)+lT(r,f)+T(r,g)+lT(r,f)+O(1)\\
                 &=2T(r,g)+2lT(r,f)+O(1)=(2l+2)T(r,g)+S(r,g)
\end{split}
\end{equation}
On the other hand, from \eqref{eq3.17X19} we have
\begin{equation}\label{eq3.23X19}
\left(\frac{f-a_2}{f-a_1}\right)^l=\frac{f_1(g-a_2)}{g-a_1} \, \, \text{and} \, \, \left(\frac{f-a_3}{f-a_1}\right)^l=\frac{f_2(g-a_3)}{g-a_1}.
\end{equation}
From the left equality of \eqref{eq3.23X19}, Lemma \ref{Nevanlinna1926}(i) and Valiron-Mohon'ko lemma (cf.\cite{Mohon'ko1971}) we have
\begin{equation}\nonumber
\begin{split}
&\quad lT(r,g)=lT(r,f)+S(r,g)=T\left(r,\left(\frac{f-a_2}{f-a_1}\right)^l\right)+S(r,g)=T\left(r,\frac{f_1(g-a_2)}{g-a_1}\right)\\
&\quad  +S(r,g)\leq T(r,f_1)+T\left(r,\frac{g-a_2}{g-a_1}\right)+S(r,g)=T(r,f_1)+T\left(r,g\right)+S(r,g),
\end{split}
\end{equation}
i.e.,
\begin{equation}\label{eq3.24X19}
(l-1)T(r,g)\leq T(r,f_1)+S(r,g).
\end{equation}
Similarly, from the right equality of \eqref{eq3.23X19}, Lemma \ref{Nevanlinna1926}(i) and the Valiron-Mohon'ko lemma (cf.\cite{Mohon'ko1971}) we have
\begin{equation}\label{eq3.25X19}
(l-1)T(r,g)\leq T(r,f_2)+S(r,g).
\end{equation}
From \eqref{eq3.24X19} and \eqref{eq3.25X19} we have
\begin{equation}\label{eq3.26X19}
2(l-1)T(r,g)\leq  T(r,f_1)+T(r,f_2)+S(r,g).
\end{equation}
From \eqref{eq3.22X19} and \eqref{eq3.26X19} we have
\begin{equation}\label{eq3.27X19}
2(l-1)T(r,g)\leq  T(r,f_1)+T(r,f_2)+S(r,g)\leq (2l+2)T(r,g)+S(r,g)
\end{equation}
Noting that $l\geq 2$ is a positive integer, we have from \eqref{eq3.21X19} and \eqref{eq3.27X19} that
\begin{equation}\label{eq3.28X19}
\overline{N}_0(r,1;f_1,f_2)\neq S(r),
\end{equation}
where and in what follows, $S(r)=S(r,g)$ is some quantity such that $S(r)=S(r,g)=o(T(r)),$ as $r\not\in E$ and $r\rightarrow\infty.$ Here $T(r)=T(r,f_1)+T(r,f_2),$ and $E\subset (0,+\infty)$ denotes a set of finite linear measure.
From \eqref{eq3.17X19}, \eqref{eq3.18X19}, \eqref{eq3.28X19} and  Lemma \ref{Zhang1999} we see that there exist two integers $s$ and $t$ with $|s|+|t|>0$ such that $f^s_1f^t_2=1.$ This together with the definitions of $f_1$ and $f_2$ in \eqref{eq3.17X19} gives
\begin{equation}\label{eq3.29X19}
\left(\frac{(g-a_1)(f-a_2)^l}{(f-a_1)^l(g-a_2)}\right)^s\left(\frac{(g-a_1)(f-a_3)^l}{(f-a_1)^l(g-a_3)}\right)^t=1.
\end{equation}
From \eqref{eq3.29X19} we have a discussion as follows:
\vskip 2mm
\par Suppose that one of $s$ and $t$ is equal to zero, say $s=0.$ Then $t\neq 0.$ Therefore, \eqref{eq3.29X19} can be rewritten into
\begin{equation}\label{eq3.30X19}
\left(\frac{f-a_3}{f-a_1}\right)^{lt}=\left(\frac{g-a_3}{g-a_1}\right)^t.
\end{equation}
 We take the Nevanlinna$'$s characteristic functions on both sides of \eqref{eq3.30X19}, and then use the Valiron-Mohon'ko lemma (cf.\cite{Mohon'ko1971}) to get
 \begin{equation}\label{eq3.31X19}
l|t|T(r,f)=|t|T(r,g)+O(1).
\end{equation}
 From \eqref{eq3.31X19}, Lemma \ref{Nevanlinna1926}(i) and the supposition that $l\geq 2$ is a positive integer, we have $T(r,f)=S(r,f).$ This is impossible.
 \vskip 2mm
\par Suppose that $s+t=0.$ Then $s=-t.$ Moreover, from $|s|+|t|> 0$ we see that $s$ and $t$ are two non-zero inters. Therefore,
\eqref{eq3.29X19} can be rewritten into
\begin{equation}\label{eq3.32X19}
\left(\frac{f-a_3}{f-a_2}\right)^{lt}=\left(\frac{g-a_3}{g-a_2}\right)^t.
\end{equation}
We take the Nevanlinna$'$s characteristic functions on both sides of \eqref{eq3.32X19} and then use the Valiron-Mohon'ko lemma (cf.\cite{Mohon'ko1971}) to get \eqref{eq3.31X19}.
From \eqref{eq3.31X19}, Lemma \ref{Nevanlinna1926}(i) and the supposition that $l\geq 2$ is a positive integer,  we have $T(r,f)=S(r,f).$ This is impossible.
 \vskip 2mm
\par Suppose that $s+t\neq 0.$ We first rewrite \eqref{eq3.29X19} into
\begin{equation}\label{eq3.33X19}
\left(\frac{(f-a_2)^s(f-a_3)^t}{(f-a_1)^{s+t}}\right)^{l}=\frac{(g-a_2)^s(g-a_3)^t}{(g-a_1)^{s+t}}.
\end{equation}
We take the Nevanlinna$'$s characteristic functions on both sides of \eqref{eq3.33X19} and then use the supposition $s+t\neq 0$ and the Valiron-Mohon'ko lemma (cf.\cite{Mohon'ko1971}) to get
\begin{equation}\label{eq3.34X19}
k|s+t|T(r,f)=|s+t|T(r,g)+O(1).
\end{equation}
From \eqref{eq3.34X19}, Lemma \ref{Nevanlinna1926}(i) and the supposition that $l\geq 2$ is a positive integer, we have $T(r,f)=S(r,f).$ This is impossible.
\vskip 2mm
\par {\bf Subcase 1.3.} Suppose that there exists some positive integer $m$ satisfying $2\leq m\leq 6,$ such that
\begin{equation}\label{eq3.35X19}
\overline{N}_{(m,1)}\left (r,a_1;f,g \right )+\overline{N}_{(m,1)}\left (r,a_2;f,g \right ) +\overline{N}_{(m,1)}\left (r,a_3;f,g \right ) \neq S(r,f).
\end{equation}
where and in what follows, $\overline{N}_{(m,1)}\left (r,a;f,g \right )$ denotes the reduced counting function of those common zeros of $f-a$ and $g-a$ in $|z|<r$ that are zeros of $f-a$ of multiplicity $m,$ and simple zeros of $g-a$ for $a\in\{a_1,a_2,a_3\}.$ Noting that \eqref{eq3.35X19} is equivalent to
\begin{equation}\label{eq3.36X19}
\overline{N}_{(1,m)}\left (r,a_1;g,f \right )+\overline{N}_{(1,m)}\left (r,a_2;g,f \right ) +\overline{N}_{(1,m)}\left (r,a_3;g,f \right ) \neq S(r,f).
\end{equation}
for the positive integer $m$ satisfying $2\leq m\leq 6.$ Next we apply \eqref{eq3.36X19} and the line of the reasoning of Subcase 1.2 to the reasoning of Subcase 1.3 with the positive integer $m\geq 2$ instead of the positive integer $l\geq 2$ and the functions $f$ and $g$ instead of $g$ and $f$ respectively, and then we get a contradiction.
\vskip 2mm
\par {\bf Case 2.} Suppose that $\Phi_{f}/\Phi_{g}$ is a non-vanishing meromorphic function in the complex plane. Without loss of generality, we suppose that
the two entire functions $\Phi_{f}$ and $\Phi_{g}$ defined as \eqref{eq1.1X1} have infinitely many zeros in the complex plane.
Next we let
\begin{equation}\label{eq3.37X19}
Z_1=\{z_{1,m}\}_{m=1}^{+\infty}\subset\Bbb{C}\setminus\{0\} \, \, \text{and} \, \, Z_2=\{z_{2,n}\}_{n=1}^{+\infty}\subset\Bbb{C}\setminus\{0\}
\end{equation}
be two infinite sequences of non-zero zeros counted with multiplicities of $\Phi_{f}$ and $\Phi_{g}$ in the complex plane respectively such that $|z_{1,m}|\rightarrow+\infty,$ as $m\rightarrow+\infty,$ and that $|z_{2,n}|\rightarrow+\infty,$ as $n\rightarrow+\infty.$ Here we mention that each point in the sequence $Z_1$ defined as \eqref{eq3.37X19} is repeated as many times as its multiplicity of a zero of $\Phi_{f},$ and that each point in the sequence $Z_2$ defined as \eqref{eq3.37X19} is repeated as many times as its multiplicity of a zero of $\Phi_{g}.$ Next we let
 \begin{equation}\label{eq3.38X19}
 \lambda_1(r)=\max\{1,N(r,Z_1)\} \, \, \text{and} \, \, \lambda_2(r)=\max\{1,N(r,Z_2)\} \, \, \text{for} \, \, r>0.
 \end{equation}
 Then, it follows from \eqref{eq3.37X19}, \eqref{eq3.38X19} and Lemma \ref{Miles1972} that
 \begin{equation}\label{eq3.39X19}
  N(r, Z_1) \leq 5\lambda_1 (4r) \, \, \text {and} \, \,  N(r, Z_2) \leq 5\lambda_2 (4r) \, \, \text{for} \, \,  r > 0,
 \end{equation}
and that
 \begin{equation}\label{eq3.40X19}
S(r,\rho;j:Z_1)=\left|\frac{1}{j}\sum\limits_{r<|z_{1,m}|\leq \rho}\left(\frac{1}{z_{1,m}}\right)^j\right|\leq \frac{A_1\lambda_1(B_1r)}{r^j}+\frac{A_1\lambda_1(B_1\rho)}{\rho^j}
\end{equation}
and
 \begin{equation}\label{eq3.41X19}
S(r,\rho;j:Z_2)=\left|\frac{1}{j}\sum\limits_{r<|z_{2,n}|\leq \rho}\left(\frac{1}{z_{2,n}}\right)^j\right|\leq \frac{A_2\lambda_1(B_2r)}{r^j}+\frac{A_2\lambda_1(B_2\rho)}{\rho^j}
\end{equation}
 for each positive integer $j$ and any pair of positive numbers $r$ and $\rho$ satisfying $0<r<\rho,$ where $A_1,$ $A_2,$ $B_1,$ $B_2$ are some four positive constants. From \eqref{eq3.38X19}-\eqref{eq3.41X19} and Definition \ref{Definition2.1} we deduce that the infinite sequence $Z_1$ defined as \eqref{eq3.37X19} has finite $\lambda_1$-density and is $\lambda_1$-balanced, and that the infinite sequence $Z_2$ defined as \eqref{eq3.37X19} has finite $\lambda_2$-density and is $\lambda_2$-balanced. Combining this with Lemma \ref{RubelTaylor1968}, we see that there exists some transcendental entire function $h_1$ of finite $\lambda_1$-type such that the infinite sequence $Z_1$ defined as \eqref{eq3.37X19} is the precise sequence of zeros counted with multiplicities of $h_1$ in the complex plane, and that there exists some transcendental entire function $h_2$ of finite $\lambda_2$-type such that the infinite sequence $Z_2$ defined as \eqref{eq3.37X19} is the precise sequence of zeros counted with multiplicities of $h_2$ in the complex plane. This implies that there exist some four positive constants, say $A_1,$ $A_2,$ $B_1,$ $B_2$ with $B_1>1$ and $B_2>1$ such that
\begin{equation}\label{eq3.42X19}
T(r,h_1)\leq A_1\lambda_1(B_1r) \, \,  \text{and}  \, \, T(r,h_2)\leq A_2\lambda_1(B_2r)  \, \, \text{for} \, \, r>0.
\end{equation}
Moreover, from \eqref{eq3.37X19} we know that the infinite sequence $Z_1=\{z_{1,m}\}_{m=1}^{+\infty}\subset\Bbb{C}\setminus\{0\}$
as mentioned above is the precise sequence of non-zero zeros counted with multiplicities of the entire function $\Phi_{f}$ and the precise sequence of zeros counted with multiplicities of the entire function $h_1$ in the complex plane, and that the infinite sequence $Z_2=\{z_{2,n}\}_{n=1}^{+\infty}\subset\Bbb{C}\setminus\{0\}$ as mentioned above is the precise sequence of non-zero zeros counted with multiplicities of the entire function $\Phi_{g}$ and the precise sequence of zeros counted with multiplicities of the entire function $h_2$ in the complex plane. This together with the supposition \eqref{eq3.2X19} and the expressions of $\Phi_{f}$ and $\Phi_{g}$ in \eqref{eq1.1X1} implies that $0$ and $\infty$ are two Picard exceptional values of $\Phi_f/\Phi_g$ and $h_2/h_1,$ and so we have that
\begin{equation}\label{eq3.43X19}
\frac{h_2\Phi_{f}}{h_1\Phi_{g}}=\frac{h_2\Phi_{f}}{h_1\Phi_{g}}=\frac {h_2f'(g-a_1)(g-a_2)(g-a_3)} {h_1g'(f-a_1)(f-a_2)(f-a_3)}
\end{equation}
 is an entire function such that $0$ and $\infty$ are its Picard  exceptional values. This reveals that $\log \frac{lh_2\Phi_{f}}{h_1\Phi_{g}}$ and $\log \frac{h_2\Phi_{f}}{mh_1\Phi_{g}}$ are also entire functions for any pair of positive integers $l$ and $m$ satisfying $1\leq l\leq 6$ and $2\leq m\leq 6$ respectively. Combining this with Lemma \ref{Conway1978p.99} and one of the properties of the continuous functions (cf.\cite[p.26,Theorem 5.8]{Conway1978}), we deduce that there exist some integers $k_l$ and $\tilde{k}_m$ that depend only upon the positive integers $l$ and $m$ satisfying $1\leq l\leq 6$ and $2\leq m\leq 6$ respectively such that
\begin{equation}\label{eq3.44X19}
\log \frac{lh_2(z)\Phi_{f}(z)}{h_1(z)\Phi_{g}(z)}=\log \frac{l\Phi_{f}(z)}{\Phi_{g}(z)}+\log \frac{h_2(z)}{h_1(z)}+2k_l\pi
\end{equation}
and
\begin{equation}\label{eq3.45X19}
 \log \frac{h_2(z)\Phi_{f}(z)}{mh_1(z)\Phi_{g}(z)}=\log \frac{\Phi_{f}(z)}{m\Phi_{g}(z)}+\log \frac{h_2(z)}{h_1(z)}+2\tilde{k}_m\pi
\end{equation}
for each $z\in \Bbb{C}$ such that $h_2(z)/h_1(z)\neq 0,\infty$ and for any pair of positive integers $l$ and $m$ satisfying $1\leq l\leq 6$ and $2\leq m\leq 6$ respectively.  Suppose that $z_{j,1,l}$ is a simple zero of $f-a_j$ and a zero of $g-a_j$ of multiplicity $l\geq 1$ for $j\in\{1,2,3\}.$ Then, by calculation we have from \eqref{eq3.7X19} that
\begin{equation}\label{eq3.46X19}
\frac{\Phi_{f}(z_{j,1,l})}{\Phi_{g}(z_{j,1,l})}=\frac {f'(z_{j,1,l})(g(z_{j,1,l})-a_1)(g(z_{j,1,l})-a_2)(g(z_{j,1,l})-a_3)} {g'(z_{j,1,l})(f(z_{j,1,l})-a_1)(f(z_{j,1,l})-a_2)(f(z_{j,1,l})-a_3)}=\frac{1}{l}.
\end{equation}
Similarly, we suppose that $z_{j,m,1}$ is a simple zero of $g-a_j$ and a zero of $f-a_j$ of multiplicity $m\geq 2$ for $j\in\{1,2,3\}.$ Then, by calculation we have from \eqref{eq3.7X19} that
\begin{equation}\label{eq3.47X19}
\frac{\Phi_{f}(z_{j,m,1})}{\Phi_{g}(z_{j,m,1})}=\frac {f'(z_{j,m,1})(g(z_{j,m,1})-a_1)(g(z_{j,m,1})-a_2)(g(z_{j,m,1})-a_3)} {g'(z_{j,m,1})(f(z_{j,m,1})-a_1)(f(z_{j,m,1})-a_2)(f(z_{j,m,1})-a_3)}=m.
\end{equation}
Without loss of generality, we suppose that $\log w$ is such a branch of the
logarithm in the open and connected region $\Bbb{C}\setminus\{0,\infty\}$ belonging the $w$-plane $\Bbb{C}$ such that $\log 1=0.$ Next we let
\begin{equation}\label{eq3.48X19}
H_{1,l}(z)=\log\frac{l\Phi_f(z)}{\Phi_g(z)} \, \, \text{and} \, \, H_{m,1}(z)=\log\frac{\Phi_f(z)}{m\Phi_g(z)}
\end{equation}
for each $z\in\Bbb{C}$ and any pair of positive integers $l$ and $m$ satisfying $1\leq l\leq 6$ and $2\leq m\leq 6$ respectively. Then, it follows from \eqref{eq3.46X19}-\eqref{eq3.48X19} that
\begin{equation}\label{eq3.49X19}
H_{1,l}(z_{j,1,l})=\log\frac{l\Phi_f(z_{j,1,l})}{\Phi_g(z_{j,1,l})}=0 \, \, \text{and} \, \, H_{m,1}(z_{j,m,1})=\log\frac{\Phi_f(z_{j,m,1})}{m\Phi_g(z_{j,m,1})}=0
\end{equation}
for each positive integer $j$ satisfying $1\leq j\leq 3$ and any pair of positive integers $l$ and $m$ satisfying $1\leq l\leq 6$ and $2\leq m\leq 6$ respectively. From \eqref{eq3.44X19}, \eqref{eq3.45X19} and \eqref{eq3.48X19} we have
\begin{equation}\label{eq3.50X19}
H_{1,l}(z)=\log \frac{lh_2(z)\Phi_{f}(z)}{h_1(z)\Phi_{g}(z)}-\log \frac{h_2(z)}{h_1(z)}-2k_l\pi
\end{equation}
and
\begin{equation}\label{eq3.51X19}
 H_{m,1}(z)=\log \frac{h_2(z)\Phi_{f}(z)}{mh_1(z)\Phi_{g}(z)}-\log \frac{h_2(z)}{h_1(z)}-2\tilde{k}_m\pi
\end{equation}
for each $z\in \Bbb{C}$ such that $h_2(z)/h_1(z)\neq 0,\infty$ and for any pair of positive integers $l$ and $m$ satisfying $1\leq l\leq 6$ and $2\leq m\leq 6$ respectively. Therefore, it follows from \eqref{eq3.5X19}, \eqref{eq3.49X19} and Lemma \ref{LiWangYi20261} that
\begin{equation}\nonumber
\begin{split}
\frac {1}{7} T(r,f)&\leq\sum\limits_{j=1}^{3}\left(\sum\limits_{l=1}^{6}\overline{N}_{(1,l)}\left(r,a_j;f,g\right)
+\sum\limits_{m=2}^{6}\overline{N}_{(m,1)}\left(r,a_j;f,g\right)\right)+S(r,f)+S(r,g)\\
&\leq 3\sum\limits_{l=1}^{6}\overline{N}\left(r,\frac{1}{H_{1,l}}\right)
+3\sum\limits_{m=2}^{6}\overline{N}\left(r,\frac{1}{H_{m,1}}\right)+S(r,f)+S(r,g)\\
&\leq 3\sum\limits_{l=1}^{6}T\left(r, \frac{1}{H_{1,l}}\right)+3\sum\limits_{m=2}^{6}T\left(r,\frac{1}{H_{m,1}}\right) +S(r,f)+S(r,g)\\
&=3\sum\limits_{l=1}^{6}m\left(r,H_{1,l}\right)+3\sum\limits_{m=2}^{6}m\left(r,H_{m,1}\right)+A_3+S(r,f)+S(r,g),
\end{split}
\end{equation}
i.e.,
\begin{equation}\label{eq3.52X19}
\frac {1}{7} T(r,f)\leq 3\sum\limits_{l=1}^{6}m\left(r,H_{1,l}\right)+3\sum\limits_{m=2}^{6}m\left(r,H_{m,1}\right) +A_3+S(r,f)+S(r,g),
\end{equation}
where and in what follows, $S(r,f)+S(r,g)$ is the error term that can be estimated as \eqref{eq3.6X19}, and $A_3$ is a real constant such that
\begin{equation}\label{eq3.53X19}
A_3=3\sum\limits_{l=1}^{6}\log\left|\left(\log\frac{l\Phi_f(0)}{\Phi_g(0)}\right)^{-1}\right|
+3\sum\limits_{m=2}^{6}\log\left|\left(\log\frac{\Phi_f(0)}{m\Phi_g(0)}\right)^{-1}\right|.
\end{equation}
From \eqref{eq3.50X19} and \eqref{eq3.51X19} we see that $m\left(r,H_{1,l}\right)$ and $m\left(r,H_{m,1}\right)$ are estimated as
\begin{equation}\label{eq3.54X19}
\begin{split}
m\left(r,H_{1,l}\right)&=m\left(r,\log \frac{lh_2\Phi_{f}}{h_1\Phi_{g}}-\log \frac{h_2}{h_1}-2k_l\pi \right)\leq m\left(r,\log \frac{h_2\Phi_{f}}{h_1\Phi_{g}}\right)\\
&\quad +m\left(r,\log\frac{h_2}{h_1}\right)+\log^{+}|2k_l\pi|+\log l+\log 3
\end{split}
\end{equation}
and
\begin{equation}\label{eq3.55X19}
\begin{split}
m\left(r,H_{m,1}\right)&=m\left(r, \log \frac{h_2\Phi_{f}}{mh_1\Phi_{g}}-\log \frac{h_2}{h_1}-2\tilde{k}_m\pi\right)\leq m\left(r,\log \frac{h_2\Phi_{f}}{h_1\Phi_{g}}\right)\\
&\quad +m\left(r,\log\frac{h_2}{h_1}\right)+\log^{+}|2\tilde{k}_m\pi|+\log 3
\end{split}
\end{equation}
for any pair of positive integers $l$ and $m$ satisfying $1\leq l\leq 6$ and $2\leq m\leq 6$ respectively. We substitute \eqref{eq3.54X19} and \eqref{eq3.55X19} into \eqref{eq3.52X19}, and then we have
\begin{equation}\label{eq3.56X19}
\begin{split}
\frac {1}{7} T(r,f)&\leq 3\sum\limits_{l=1}^{6}m\left(r,H_{1,l}\right)+3\sum\limits_{m=2}^{6}m\left(r,H_{m,1}\right) +A_3+S(r,f)+S(r,g)\\
&\leq 33m\left(r,\log \frac{h_2\Phi_{f}}{h_1\Phi_{g}}\right)+33m\left(r,\log \frac{h_2}{h_1}\right)+A_4+S(r,f)+S(r,g),\\
\end{split}
\end{equation}
where $A_4$ is a real constant such that
\begin{equation}\nonumber
A_4=A_3+3\sum\limits_{l=1}^6\log^{+}|2k_l\pi|+3\sum\limits_{m=2}^6\log^{+}|2\tilde{k}_m\pi|+3\sum\limits_{l=1}^6\log l+33\log 3
\end{equation}
with the real constant $A_3$ defined as \eqref{eq3.53X19}. Next we estimate $m\left(r,\log \frac{h_2\Phi_{f}}{h_1\Phi_{g}}\right)$ and
$m\left(r,\log \frac{h_2}{h_1}\right)$ in \eqref{eq3.56X19}. Firstly we estimate $m\left(r,\log \frac{h_2\Phi_{f}}{h_1\Phi_{g}}\right)$ as follows:
on the one hand, from the known result that $0$ and $\infty$ are two Picard exceptional values of the entire function $\frac{h_2\Phi_{f}}{h_1\Phi_{g}}$ in \eqref{eq3.43X19}, we see that $\log \frac{h_2\Phi_{f}}{h_1\Phi_{g}}$ is also an entire function. Combining this with \eqref{eq3.1X19}, \eqref{eq3.2X19}, \eqref{eq3.43X19} and Lemma \ref{Chuang1951}, we have
\begin{equation}\label{eq3.57X19}
\begin{split}
&\quad m\left(r,\log \frac{h_2\Phi_{f}}{h_1\Phi_{g}}\right)=T\left(r,\log \frac{h_2\Phi_f}{h_1\Phi_{g}}\right)\\
&<A_5T\left(\tau r,\frac{h'_2}{h_2}+\frac{\Phi'_f}{\Phi_f}-\frac{h'_1}{h_1}-\frac{\Phi'_g}{\Phi_g}\right)+\log^{+}(\tau r)+4+\log^{+}\left|\log \frac{h_2(0)\Phi_f(0)}{h_1(0)\Phi_g(0)}\right|\\
&=A_5 m\left(\tau r,\frac{h'_2}{h_2}+\frac{f''}{f'}+\sum\limits_{j=1}^3\frac{g'}{g-a_j}-\frac{h'_1}{h_1}-\frac{g''}{g'}-\sum\limits_{j=1}^{3}\frac{f'}{f-a_j}\right)
+\log^{+}(\tau r)+4\\
&\quad +\log^{+}\left|\log \frac{h_2(0)\Phi_{f}(0)}{h_1(0)\Phi_g(0)}\right|\leq  A_6\left(m\left(\tau r,\frac{f''}{f'}\right)+m\left(\tau r,\frac{g''}{g'}\right)\right)\\
 \end{split}
 \end{equation}
 \begin{equation}\nonumber
 \begin{split}
&\quad+A_6\sum\limits_{j=1}^{3}\left(m\left(\tau r,\frac{f'}{f-a_j}\right)+m\left(\tau r,\frac{g'}{g-a_j}\right)\right)+A_6m\left(\tau r,\frac{h'_1}{h_1}\right)\\
&\quad+A_6m\left(\tau r,\frac{h'_2}{h_2}\right)+\log^{+}r+A_6\\
\end{split}
\end{equation}
for a given positive constant $\tau>1$ and any $r>0.$ Here $A_5$ is a positive constant depending only upon the positive constant $\tau,$ and $A_6$ is some large positive constant such that $A_6>A_5$. On the one hand, from \eqref{eq2.11X19} in Remark \ref{remark2.4} we have
\begin{equation}\label{eq3.58X19}
\begin{split}
&\quad m\left(\tau r,\frac{f''}{f'}\right)+m\left(\tau r,\frac{g''}{g'}\right)+\sum\limits_{j=1}^{3}\left(m\left(\tau r,\frac{f'}{f-a_j}\right)+m\left(\tau r,\frac{g'}{g-a_j}\right)\right)\\
&\leq C_3\left(\log^{+}T(\rho,f)+\log^{+}\log^{+}T(\rho,f)+\log^{+}\frac{\rho}{\rho-\tau r}+1\right)\\
&=C_3\left(\log^{+}T(\tau^2 r,f)+\log^{+}\log^{+}T(\tau^2 r,f)+\log^{+}\frac{\tau^2 r}{\tau^2 r-\tau r}+1\right)\\
&=C_3\left(\log^{+}T(\tau^2 r,f)+\log^{+}\log^{+}T(\tau^2 r,f)+\log^{+}\frac{\tau }{\tau -1 }+1\right)
\end{split}
\end{equation}
for the positive numbers $r$ and $\rho$ such that $1<r_0<r<\tau r<\tau^2 r<R<+\infty$ with $\rho=\tau^2 r$ and the positive constant $\tau>1.$
On the other hand, from the definitions of the infinite sequences $Z_1$ and $Z_2$ of \eqref{eq3.37X19} of non-zero zeros of the entire functions $\Phi_f$  and $\Phi_g$ in \eqref{eq1.1X1} respectively, we have from \eqref{eq3.38X19}-\eqref{eq3.41X19} that there exists some large positive number, say $r_0>1$ such that
\begin{equation}\label{eq3.59X19}
 \lambda_1(r)=\max\{1,N(r,Z_1)\}=N(r,Z_1) \, \, \text{and} \, \,  N(r, Z_1) \leq  5\lambda_1(4r)
 \end{equation}
and
\begin{equation}\label{eq3.60X19}
 \lambda_2(r)=\max\{1,N(r,Z_2)\}=N(r,Z_2) \, \, \text{and} \, \,  N(r, Z_2) \leq  5\lambda_5(4r)
 \end{equation}
 for $r> r_0>1.$ Noting that the infinite sequence $Z_1$ in \eqref{eq3.37X19} counted with multiplicities is the non-zero zeros of the entire function $\Phi_f$ in the complex plane, and the precise sequence of zeros counted with multiplicities of the entire function $h_1$ in the complex plane, while the infinite sequence $Z_2$ in \eqref{eq3.37X19} counted with multiplicities is the non-zero zeros of the entire function $\Phi_g$ in the complex plane, and the precise sequence of zeros counted with multiplicities of the entire function $h_2$ in the complex plane, we have from \eqref{eq3.4X19}, the left formula of \eqref{eq3.42X19} and the left formulas of \eqref{eq3.59X19} and \eqref{eq3.60X19} that
\begin{equation}\label{eq3.61X19}
\begin{split}
&\quad T(r,h_1)\leq A_1\lambda_1(B_1r)= A_1N\left(B_1r,\frac{1}{h_1}\right)\leq A_1N\left(B_1r,\frac{1}{\Phi_f}\right)\leq A_1N\left(B_1r,\frac{1}{\Psi}\right)\\
&=A_1 C_3\left(\log^{+}T(B^2_1r,f)+\log^{+}\log^{+}T(B^2_1r,f)+\log^{+}\frac{B^2_1r}{B^2_1r-B_1r}+1\right)\\
&<A_7 \left(\log^{+}T(B^2_1r,f)+\log^{+}\log^{+}T(B^2_1r,f)+1\right)\\
\end{split}
\end{equation}
for the positive number $r$ such that $1<r_0<r<B_1r<B^2_1r<R<+\infty$ with the positive constant $B_1>1$ as mentioned in the line before \eqref{eq3.42X19}.
Here and in what follows, $A_7$ is a large positive constant depending only upon the positive constants $A_1,$ $C_3$ and the constant $B_1>1$ such that $A_7>A_1C_3.$ Similarly, from  \eqref{eq3.4X19}, the left formula of \eqref{eq3.42X19} and the left formulas of \eqref{eq3.59X19} and \eqref{eq3.60X19} we have
\begin{equation}\label{eq3.62X19}
\begin{split}
T(r,h_2)&\leq A_2\lambda_2(B_2r)= A_2N\left(B_2r,\frac{1}{h_2}\right)\leq A_2N\left(B_2r,\frac{1}{\Phi_g}\right)\leq A_2N\left(B_2r,\frac{1}{\Psi}\right)\\
&\leq  A_2 C_3\left(\log^{+}T(B^2_2r,f)+\log^{+}\log^{+}T(B^2_2r,f)+\log^{+}\frac{B^2_2r}{B^2_2r-B_2r}+1\right)\\
&<A_8 \left(\log^{+}T(B^2_2r,f)+\log^{+}\log^{+}T(B^2_2r,f)+1\right)\\
\end{split}
\end{equation}
for the positive number $r$ such that $1<r_0<r<B_2r<B^2_2r<R<+\infty.$ Here and in what follows, $A_8$ is a large positive constant depending only upon the positive constants $A_2,$ $C_3$ and the positive constant $B_2>1$ as mentioned in the line before \eqref{eq3.42X19}, such that $A_8>A_2C_3.$ From \eqref{eq3.61X19} and Lemma \ref{Heittokangas2004} we see that there exist some two large positive constants, say $A_{9}$ and $r_0$ with $A_{9}>1$ and $r_0>1$ such that
\begin{equation}\label{eq3.63X19}
\begin{split}
&\quad m\left(\tau r,\frac{h'_1}{h_1}\right)\leq A_9\log^{+}\left(\frac{T(\rho,h_1)}{\tau r}\cdot\frac{\rho}{\rho-\tau r}\right)\\
&= A_9\log^{+}\left(\frac{T(B_1\tau r,h_1)}{\tau r}\cdot\frac{B_1\tau r}{B_1\tau r-\tau r}\right)\\
&\leq A_{9}\log^{+}\left(\frac{A_7B_1\left(\log^{+}T(B^3_1\tau r,f)+\log^{+}\log^{+}T(B^3_1\tau r,f)+1\right)}{(B_1-1)\tau r}\right)\\
&\leq A_{10}\left(\log^{+}\log^{+}T(B^3_1\tau r,f)+\log^{+}\log^{+}\log^{+}T(B^3_1\tau r,f)+1\right)
\end{split}
\end{equation}
for the positive numbers $r$ and $\rho$ with $\rho=B_1\tau r$ such that
 \begin{equation}\label{eq3.64X19}
 1<r_0<r<\tau r<B_1\tau r<B^2_1\tau r<B^3_1\tau r<R<+\infty.
   \end{equation}
 Here and in what follows, $A_{10}$ is some positive constant depending only upon the positive constants $A_7,$ $A_9,$ $B_1$ and $\tau$ with $A_9>1,$ $B_1>1$ and $\tau>1$ such that $A_{10}\geq A_9.$ Similarly, from  \eqref{eq3.62X19} and Lemma \ref{Heittokangas2004} we see that there exist some two large positive constants, say $A_{11}$ and $r_0$ with $A_{11}>1$ and $r_0>1$ such that
 \begin{equation}\label{eq3.65X19}
 \begin{split}
m\left(\tau r,\frac{h'_2}{h_2}\right)\leq A_{11}\left(\log^{+}\log^{+}T(B^3_2\tau r,f)+\log^{+}\log^{+}\log^{+}T(B^3_2\tau r,f)+1\right)
\end{split}
\end{equation}
for the positive number $r$ such that
 \begin{equation}\label{eq3.66X19}
 1<r_0<r<\tau r<B_2\tau r<B^2_2\tau r<B^3_2\tau r<R<+\infty
   \end{equation}
 Here and in what follows, $A_{11}$ is some large positive constant depending only upon the constants $A_8,$ $B_2$ and $\tau$ with $B_2>1$ and $\tau>1.$  From \eqref{eq3.63X19}-\eqref{eq3.66X19} we have
 \begin{equation}\label{eq3.67X19}
\begin{split}
&\quad m\left(\tau r,\frac{h'_1}{h_1}\right)+ m\left(\tau r,\frac{h'_2}{h_2}\right)\\
&\leq A_{10}\left(\log^{+}\log^{+}T(B^3_1\tau r,f)+\log^{+}\log^{+}\log^{+}T(B^3_1\tau r,f)+1\right)\\
&\quad + A_{11}\left(\log^{+}\log^{+}T(B^3_2\tau r,f)+\log^{+}\log^{+}\log^{+}T(B^3_2\tau r,f)+1\right)\\
\end{split}
\end{equation}
 for the positive number $r$ such that \eqref{eq3.64X19} and \eqref{eq3.66X19} hold. Noting that $f$ is a non-constant meromorphic function, we have that $T(r,f)$ is a  monotonically increasing, continuous and unbounded function for $0 < r < +\infty$ (cf.\cite[p.8, Corollary 1]{Hayman1964}). Without loss of generality, next we suppose that the positive constants $B_1$ and $B_2$ in \eqref{eq3.42X19} and the positive constant $\tau$ in \eqref{eq3.57X19} satisfy  $1<\tau<\max\{B_1, B_2\}=:B_3.$ Then we substitute \eqref{eq3.58X19} and \eqref{eq3.67X19} into \eqref{eq3.57X19}, and then we have
\begin{equation}\label{eq3.68X19}
\begin{split}
&\quad m\left(r,\log \frac{h_2\Phi_{f}}{h_1\Phi_{g}}\right)< C_3A_6\left(\log T(\tau^2 r,f)+\log\log T(\tau^2 r,f)+\log\frac{\tau }{\tau -1 }+1\right)\\
&\quad + A_6A_{10}\left(\log\log T(B^3_1\tau r,f)+\log\log\log T(B^3_1\tau r,f)+1\right)\\
&\quad + A_6A_{11}\left(\log\log T(B^3_2\tau r,f)+\log\log\log T(B^3_2\tau r,f)+1\right)\\
&\quad+\log r+A_6\leq A_{12}\log T(B^3_3\tau r,f)+\log r+A_{12}
\end{split}
\end{equation}
for the positive number $r$ such such that \eqref{eq3.64X19} and \eqref{eq3.66X19} hold. Here and in what follows, $A_{12}$ is a large positive constant depending only upon the positive constants $C_3,$ $A_6,$ $A_{10},$ $A_{11}$ and $\tau.$ Secondly we estimate $m\left(r,\log \frac{h_2}{h_1}\right)$ on the right hand side of \eqref{eq3.56X19}: for each complex number $z\in\Bbb{C}$ with $|z|=r$ satisfying \eqref{eq3.64X19} and \eqref{eq3.66X19} such that
\begin{equation}\label{eq3.69X19}
 \log\frac{h_2(z)}{h_1(z)}\not\in\{0,\infty\} \, \, \text{and} \, \, \log \frac{h_2(z)}{h_1(z)}=\log\left|\frac{h_2(z)}{h_1(z)}\right|+i\arg \frac{h_2(z)}{h_1(z)}
\end{equation}
From the third equality of \eqref{eq3.69X19} we have
\begin{equation}\label{eq3.70X22}
\left|\log \frac{h_2(z)}{h_1(z)}\right|\leq \left|\log\left|\frac{h_2(z)}{h_1(z)}\right|\right|+\left|\arg \frac{h_2(z)}{h_1(z)}\right|=\log^{+}\left|\frac{h_2(z)}{h_1(z)}\right|+\log^{+}\left|\frac{h_1(z)}{h_2(z)}\right|+\left|\arg \frac{h_2(z)}{h_1(z)}\right|
\end{equation}
for each complex number $z\in\Bbb{C}$ satisfying the left inequality of \eqref{eq3.69X19} with $|z|=r$ satisfying \eqref{eq3.64X19} and \eqref{eq3.66X19}.
From Cherry-Ye\cite[p.16, Proposition 1.2.2]{Cherry2001} we have
\begin{equation}\label{eq3.71X23}
\log ^{+}x=\log\max\{1, x\}=\log \frac{x+1+|x-1|}{2}\leq \frac{1}{2}\log (1+x^2)  \, \, \text{for} \, \, x\in [0,+\infty).
\end{equation}
From \eqref{eq3.70X22} and \eqref{eq3.71X23} we have
\begin{equation}\label{eq3.72X25}
\begin{split}
&\quad\log^{+}\left|\log \frac{h_2(z)}{h_1(z)}\right|\leq \log^{+}\log^{+}\left|\frac{h_2(z)}{h_1(z)}\right|+\log^{+}\log^{+}\left|\frac{h_1(z)}{h_2(z)}\right|+\log^{+}\left|\arg \frac{h_2(z)}{h_1(z)}\right|\\
&\leq \log^{+}\log^{+}\left|\frac{h_2(z)}{h_1(z)}\right|+\log^{+}\log^{+}\left|\frac{h_1(z)}{h_2(z)}\right|+\frac{1}{2}\log \left(1+\left(\arg \frac{h_2(z)}{h_1(z)}\right)^2\right) +\log 3\\
&\leq  \log^{+}\log^{+}\left|\frac{h_2(z)}{h_1(z)}\right|+\log^{+}\log^{+}\left|\frac{h_1(z)}{h_2(z)}\right|+\frac{1}{2}\log \left(1+\left(\arg \frac{h_2(z)}{h_1(z)}\right)^2\right)+\log 3\\
&\leq \log^{+}\left(\log^{+}\left|h_2(z)\right|+\log^{+}\left|\frac{1}{h_1(z)}\right|\right)+\log^{+}\left(\log^{+}\left|h_1(z)\right|+\log^{+}\left|\frac{1}{h_2(z)}\right|\right)\\
&\quad +\frac{1}{2}\log \left(1+\left(\arg \frac{h_2(z)}{h_1(z)}\right)^2\right)+\log 3\leq \log^{+}\log^{+}\left|h_2(z)\right|
+\log^{+}\log^{+}\left|\frac{1}{h_1(z)}\right|\\
&\quad+\log^{+}\log^{+}\left|h_1(z)\right|+\log^{+}\log^{+}\left|\frac{1}{h_2(z)}\right|+\frac{1}{2}\log \left(1+\left(\arg \frac{h_2(z)}{h_1(z)}\right)^2\right)+\log 12
\end{split}
\end{equation}
for each complex number $z\in\Bbb{C}$ satisfying the left inequality of \eqref{eq3.69X19} with $|z|=r$ satisfying \eqref{eq3.64X19} and \eqref{eq3.66X19}.
On the other hand, from Lemma \ref{Yang1993} we see that there exists some positive constant $A_{13}$ satisfying $A_{13}>e$ such that
\begin{equation}\label{eq3.73X25}
A_{13}\log^{+}\log^{+}\left|\frac{1}{h_1(z)}\right|\leq\log^{+}|h_1(z)| +\log^{+}\frac{1}{|h_1(z)|}+A_{13}\left(\log A_{13}-1\right)
\end{equation}
for the complex number $z\in \Bbb{C}$ satisfying $0<|h_1(z)|<1$ and the left inequality of \eqref{eq3.69X19} with $|z|=r$ satisfying \eqref{eq3.64X19} and \eqref{eq3.66X19}. Noting that
\begin{equation}\label{eq3.74X25}
A_{13}\log^{+}\log^{+}\left|\frac{1}{h_1(z)}\right|\leq 0
\end{equation}
for the complex number $z\in\Bbb{C}$ satisfying $1\leq |h_1(z)|<+\infty$ and the left inequality of \eqref{eq3.69X19} with $|z|=r$ satisfying \eqref{eq3.64X19} and \eqref{eq3.66X19}, we have from \eqref{eq3.73X25} and \eqref{eq3.74X25} that
\begin{equation}\label{eq3.75X25}
A_{13}\log^{+}\log^{+}\left|\frac{1}{h_1(z)}\right|\leq\log^{+}|h_1(z)| +\log^{+}\frac{1}{|h_1(z)|}+A_{13}\left(\log A_{13}-1\right)
\end{equation}
and
\begin{equation}\label{eq3.76X25}
A_{13}\log^{+}\log^{+}\left|h_1(z)\right|\leq\log^{+}|h_1(z)| +\log^{+}\frac{1}{|h_1(z)|}+A_{13}\left(\log A_{13}-1\right)
\end{equation}
for the complex number $z\in\Bbb{C}$ satisfying $0<|h_1(z)|<+\infty$ and the left inequality of \eqref{eq3.69X19} with $|z|=r$ satisfying \eqref{eq3.64X19} and \eqref{eq3.66X19}. Similarly, from Lemma \ref{Yang1993} we see that there exists some positive constant $A_{14}$ satisfying $A_{14}>e$ such that
\begin{equation}\label{eq3.77X25}
A_{14}\log^{+}\log^{+}\left|\frac{1}{h_2(z)}\right|\leq\log^{+}|h_2(z)| +\log^{+}\frac{1}{|h_2(z)|}+A_{14}\left(\log A_{14}-1\right)
\end{equation}
and
\begin{equation}\label{eq3.78X25}
A_{14}\log^{+}\log^{+}\left|h_2(z)\right|\leq\log^{+}|h_2(z)| +\log^{+}\frac{1}{|h_2(z)|}+A_{14}\left(\log A_{14}-1\right)
\end{equation}
for the complex number $z\in\Bbb{C}$ satisfying $0<|h_2(z)|<+\infty$ and the left inequality of \eqref{eq3.69X19} with $|z|=r$ satisfying \eqref{eq3.64X19} and \eqref{eq3.66X19}. From \eqref{eq3.72X25} and \eqref{eq3.75X25}-\eqref{eq3.78X25} we have
\begin{equation}\label{eq3.79X22}
\begin{split}
&\quad\log^{+}\left|\log\frac{h_2(z)}{h_1(z)}\right|\leq \sum\limits_{l=1}^2\frac{2}{A_{12+l}}\left(\log^{+}|h_l(z)|+\log^{+}\frac{1}{|h_l(z)|}\right)\\
&\quad +\frac{1}{2}\log \left(1+\left(\arg \frac{h_2(z)}{h_1(z)}\right)^2\right)+2\sum\limits_{l=1}^2\left(\log A_{12+l}-1\right)+\log 12\\
 \end{split}
\end{equation}
for each complex number $z\in\Bbb{C}$ satisfying $0<|h_l(z)|<+\infty$ and the left inequality of \eqref{eq3.69X19} with $l\in\{1,2\}$
and $|z|=r$ satisfying \eqref{eq3.64X19} and \eqref{eq3.66X19}. Whether one of $h_1,$ $h_2$ and $\log\frac{h_2}{h_1}$ has zeros or not on the circle $|z|=r$
for the positive number $r$ satisfying \eqref{eq3.64X19} and \eqref{eq3.66X19}, we have from \eqref{eq3.79X22} and the similar reasoning of \cite[pp.2-3, Case 2 of the proof of Theorem 1.1]{Hayman1964} that
\begin{equation}\label{eq3.80X22}
\begin{split}
 &\quad   m\left(r, \log\frac{h_2}{h_1}\right)\leq A_{15}\sum\limits_{l=1}^2\left(m\left(r,h_1\right)+m\left(r,\frac{1}{h_l}\right)\right)\\
 &\quad +\frac{\int_{0}^{2\pi}\log \left(1+\left(\arg \frac{h_2(re^{i\theta})}{h_1(re^{i\theta})}\right)^2\right)d\theta}{4\pi}+A_{15}
  \end{split}
\end{equation}
for each complex number
\begin{equation}\label{eq3.81X25}
z=re^{i\theta}\in D=:\left\{z\in\Bbb{C}:\frac{h_2(z)}{h_1(z)}\neq 0 \, \, \text{and} \, \, \frac{h_2(z)}{h_1(z)}\neq \infty \right\} \, \, \text{with} \, \, r>0.
\end{equation}
Here and in what follows $A_{15}$ is a positive constant such that
\begin{equation}\nonumber
A_{15}>\max\left\{1,\frac{2}{A_{13}}, \frac{2}{A_{14}}, 2\sum\limits_{l=1}^2\left(\log A_{12+l}-1\right)-1+\log 12\right\}.
\end{equation}
Whether one of $h_1/h_2,$ $h_2/h_1$  has zeros or not on the circle $|z|=r$ for $r>0,$ we use the similar reasoning of \cite[pp.2-3, Case 2 of the proof of Theorem 1.1]{Hayman1964} to deduce that the definite integral in the left hand side of the equality
\begin{equation}\label{eq3.82X24}
\begin{split}
&\quad \frac{1}{4\pi}\int_{0}^{2\pi}\log\left(1+\left(\arg \frac{h_2(re^{i\theta})}{h_1(re^{i\theta})}\right)^2\right)d\theta
=\frac{1}{4\pi i}\int_{\gamma_{0,r}}\frac{1}{z}\log\left(1+\left(\arg \frac{h_2(z)}{h_1(z)}\right)^2\right)dz
\end{split}
\end{equation}
is convergent to a real number depending only upon the given positive number $r,$ where and in what follows
\begin{equation}\nonumber
\gamma_{0,r}=\left\{\zeta\in\Bbb{C}:\zeta=re^{i\theta} \, \text{ with the real variable}  \, \, \theta \, \, \text{from}  \, \, 0  \, \text{ to} \, 2\pi\right \}
\end{equation}
denotes a circle with positive orientation such that its radius is equal to $r$ and its center is the origin point. Whether one of $h_1/h_2,$ $h_2/h_1$  has zeros or not on the circle $|z|=r$ for $r>0,$ we can find a sufficiently small positive number, say $\varepsilon$ such that $h_1/h_2$ and  $h_2/h_1$ have no zeros on the circle $|z|=r+\varepsilon,$ and so $\{z\in\Bbb{C}:|z|=r+\varepsilon\}\subset D$
with the region $D$ defined as \eqref{eq3.81X25}, such that we have from \eqref{eq3.82X24} and the continuity that
\begin{equation}\label{eq3.83X25}
\left|\frac{\int_{0}^{2\pi}\log\left(1+\left(\arg \frac{h_2(re^{i\theta})}{h_1((re^{i\theta})}\right)^2\right)d\theta}{4\pi}-\frac{\int_{0}^{2\pi}\log\left(1+\left(\arg \frac{h_2((r+\varepsilon)e^{i\theta})}{h_1((r+\varepsilon)e^{i\theta})}\right)^2\right)d\theta}{4\pi}\right|<1,
\end{equation}
and such that we have from the argument principle (cf.\cite[p.152,Theorem 18]{Ahlfors1979}) that
\begin{equation}\label{eq3.84X25}
\begin{split}
&\quad\Delta_{\gamma_{0,r}}\arg \frac{ h_2}{h_1}=\frac{1}{2\pi i}\int_{\gamma_{0,r}}d\left(\log \frac{h_2(z)}{h_1(z)}\right)=\frac{1}{2\pi i}\int_{\gamma_{0,r+\varepsilon}}d\left(\log \frac{h_2(z)}{h_1(z)}\right)\\
&=\Delta_{\gamma_{0,r+\varepsilon}}\arg \frac{ h_2}{h_1}=n\left(r,\frac{h_1}{h_2}\right)-n\left(r,\frac{h_2}{h_1}\right).
\end{split}
\end{equation}
Next we denote by
\begin{equation}\nonumber
\begin{split}
&\quad\gamma_{\theta_{r+\varepsilon}-2\pi,r+\varepsilon}\\
&=\left\{ z\in\Bbb{C}:z = (r+\varepsilon)e^{i\theta} \, \ \text {with the real variable} \, \,  \theta \, \, \text{from} \, \,  \theta_{r+\varepsilon}-2\pi \, \, \text{to} \, \, \theta_{r+\varepsilon}\right\}
\end{split}
\end{equation}
 a circle with the positive orientation such that its radius is equal to the positive number $r+\varepsilon$ and its center is the origin point. Then it follows \eqref{eq3.84X25} that
\begin{equation}\label{eq3.85X25}
\begin{split}
&\quad \Delta_{\gamma_{0,r}}\arg \frac{ h_2}{h_1}=\Delta_{\gamma_{0,r+\varepsilon}}\arg \frac{ h_2}{h_1}=\Delta_{\gamma_{\theta_{r+\varepsilon}-2\pi,r+\varepsilon}}\arg \frac{h_2}{h_1}\\
&=\frac{1}{2\pi i}\int_{\gamma_{\theta_{r+\varepsilon}-2\pi,r+\varepsilon}}d\left(\log \frac{h_2(z)}{h_1(z)}\right)=\arg \frac{h_2((r+\varepsilon)e^{i\theta_{r+\varepsilon}})}{h_1((r+\varepsilon)e^{i\theta_{r+\varepsilon}})}\\
&\quad-\arg \frac{h_2((r+\varepsilon)e^{i(\theta_{r+\varepsilon}-2\pi)})}{h_1((r+\varepsilon)e^{i(\theta_{r+\varepsilon}-2\pi)})}=n\left(r,\frac{h_1}{h_2}\right)-n\left(r,\frac{h_2}{h_1}\right)
\end{split}
\end{equation}
with $\arg \frac{h_2((r+\varepsilon)e^{i\theta_{r+\varepsilon}})}{h_1((r+\varepsilon)e^{i\theta_{r+\varepsilon}})}\in [-\pi,\pi)$ and $\arg \frac{h_2((r+\varepsilon)e^{i(\theta_{r+\varepsilon}-2\pi)})}{h_1((r+\varepsilon)e^{i(\theta_{r+\varepsilon}-2\pi)})}\in [2k_0\pi-\pi,2k_0\pi+\pi)$ with $k_0$ being some fixed integer, while $n\left(r,h_1/h_2\right)$ denotes the numbers of zeros counted with multiplicities of $h_2/h_1$ in the open disc $|z|<r,$ and $n\left(r,h_2/h_1\right)$ denotes the numbers of poles counted with multiplicities of $h_2/h_1$ in the open disc $|z|<r.$
From the known result that $h_1/h_2$ and $h_2/h_1$ have no zeros on the circle $|z|=r+\varepsilon,$ we deduce that $\log \left(1+\left(\arg \frac{h_2((r+\varepsilon)e^{i\theta})}{h_1((r+\varepsilon)e^{i\theta})}\right)^2\right)$ is a continuous real-valued function of the real variable $\theta$ defined on the interval $[0.2\pi]$ for the positive number $r+\varepsilon.$ Combining this with the mean value theorem for integrals, we have that there exists some real number $\theta_{r+\varepsilon}\in (0,2\pi)$ such that
\begin{equation}\label{eq3.86X25}
\begin{split}
&\quad \frac{\int_{0}^{2\pi}\log\left(1+\left(\arg \frac{h_2((r+\varepsilon)e^{i\theta})}{h_1((r+\varepsilon)e^{i\theta})}\right)^2\right)d\theta}{4\pi}
=\frac{\log\left(1+\left(\arg \frac{h_2((r+\varepsilon)e^{i\theta_{r+\varepsilon}})}{h_1((r+\varepsilon)e^{i\theta_{r+\varepsilon}})}\right)^2\right)}{2}
\end{split}
\end{equation}
From \eqref{eq3.85X25}, \eqref{eq3.86X25}, the argument principle (cf.\cite[p.152,Theorem 18]{Ahlfors1979}) and one of the properties of the continuous functions (cf.\cite[p.26,Theorem 5.8]{Conway1978}) we have
\begin{equation}\label{eq3.87X24}
\begin{split}
&\quad \frac{\int_{0}^{2\pi}\log\left(1+\left(\arg \frac{h_2((r+\varepsilon)e^{i\theta})}{h_1((r+\varepsilon)e^{i\theta})}\right)^2\right)d\theta}{4\pi}
=\frac{\log\left(1+\left(\arg \frac{h_2((r+\varepsilon)e^{i\theta_{r+\varepsilon}})}{h_1((r+\varepsilon)e^{i\theta_{r+\varepsilon}})}\right)^2\right)}{2}\\
&\leq \frac{1}{2}\log\left(1+\left(\left|\Delta_{\gamma_{\theta_{r+\varepsilon}-2\pi,r+\varepsilon}}\arg \frac{h_2}{h_1}\right|+\left|\arg \frac{h_2((r+\varepsilon)e^{i(\theta_{r+\varepsilon}-2\pi)})}{h_1((r+\varepsilon)e^{i(\theta_{r+\varepsilon}-2\pi)})}\right|\right)^2\right)\\
&=\frac{1}{2}\log\left(1+\left(\left|n\left(r,\frac{h_1}{h_2}\right)-n\left(r,\frac{h_2}{h_1}\right)\right|+\left|\arg \frac{h_2((r+\varepsilon)e^{i(\theta_{r+\varepsilon}-2\pi )})}{h_1((r+\varepsilon)e^{i(\theta_{r+\varepsilon}-2\pi)})}\right|\right)^2\right)\\
&\leq \frac{1}{2}\log\left(1+\left(n\left(r,\frac{h_1}{h_2}\right)+n\left(r,\frac{h_2}{h_1}\right)+(2|k_0|+1)\pi\right)^2\right)\\
&\leq \log \left(n\left(r,\frac{h_1}{h_2}\right)+n\left(r,\frac{h_2}{h_1}\right)+(2|k_0|+1)\pi\right)+\frac{\log 2}{2}\\
&\leq \frac{1}{(2|k_0|+1)\pi}\left(n\left(r,\frac{h_1}{h_2}\right)+n\left(r,\frac{h_2}{h_1}\right)\right)+\log \left((2|k_0|+1)\pi\right)+\frac{\log 2}{2}
\end{split}
\end{equation}
for $r>0.$ Whether one of $h_1/h_2,$ $h_2/h_1$  has zeros or not on the circle $|z|=r$ for $r>0,$ we have from \eqref{eq3.83X25} and \eqref{eq3.87X24} that
\begin{equation}\label{eq3.88X24}
\begin{split}
&\quad \frac{\int_{0}^{2\pi}\log\left(1+\left(\arg \frac{h_2(re^{i\theta})}{h_1((re^{i\theta})}\right)^2\right)d\theta}{4\pi}<1+ \frac{\int_{0}^{2\pi}\log\left(1+\left(\arg \frac{h_2((r+\varepsilon)e^{i\theta})}{h_1((r+\varepsilon)e^{i\theta})}\right)^2\right)d\theta}{4\pi}\\
&\leq \frac{1}{(2|k_0|+1)\pi}\left(n\left(r,\frac{h_1}{h_2}\right)+n\left(r,\frac{h_2}{h_1}\right)\right)+\log \left((2|k_0|+1)\pi\right)+\frac{\log 2}{2}+1
\end{split}
\end{equation}
for $r>0.$ From Hayman \cite[p.37]{Hayman1964} we have
\begin{equation}\label{eq3.89X24}
n\left(r,\frac{h_1}{h_2}\right)+n\left(r,\frac{h_2}{h_1}\right)\leq \frac{2\tau r }{\tau r -r}\left(T\left(\tau r,\frac{h_2}{h_1}\right)+\log^{+}\left|\frac{h_1(0)}{h_2(0)}\right|\right)
\end{equation}
for $r>0.$ Here and in what follows, $\tau>1$ is a positive constant defined as \eqref{eq3.64X19}. We substitute \eqref{eq3.89X24} into \eqref{eq3.88X24}, and then we have
\begin{equation}\label{eq3.90X24}
\begin{split}
&\quad \frac{\int_{0}^{2\pi}\log\left(1+\left(\arg \frac{h_2(re^{i\theta})}{h_1((re^{i\theta})}\right)^2\right)d\theta}{4\pi}<1+ \frac{\int_{0}^{2\pi}\log\left(1+\left(\arg \frac{h_2((r+\varepsilon)e^{i\theta})}{h_1((r+\varepsilon)e^{i\theta})}\right)^2\right)d\theta}{4\pi}\\
&\leq \frac{2\tau \left(T\left(\tau r,\frac{h_2}{h_1}\right)+\log^{+}\left|\frac{h_1(0)}{h_2(0)}\right|\right)}{(2|k_0|+1)(\tau-1)\pi}+\log \left((2|k_0|+1)\pi\right)+\frac{\log 2}{2}+1
\end{split}
\end{equation}
for $r>0.$ From \eqref{eq3.80X22} and \eqref{eq3.90X24} we have that there exists a large positive constant, say $r_0>1$ such that
\begin{equation}\nonumber
\begin{split}
 &\quad   m\left(r, \log\frac{h_2}{h_1}\right)\leq A_{15}\sum\limits_{l=1}^2\left(m\left(r,h_1\right)+m\left(r,\frac{1}{h_l}\right)\right)\\
 &\quad +\frac{2\tau \left(T\left(\tau r,\frac{h_2}{h_1}\right)+\log^{+}\left|\frac{h_1(0)}{h_2(0)}\right|\right)}{(2|k_0|+1)(\tau-1)\pi}+A_{15}+\log \left((2|k_0|+1)\pi\right)+\frac{\log 2}{2}+1\\
 &\leq A_{15}\sum\limits_{l=1}^2\left(T\left(r,h_1\right)+T\left(r,\frac{1}{h_l}\right)\right)+\frac{2\tau \left(T\left(\tau r,h_2\right)+T\left(\tau r,\frac{1}{h_1}\right)+\log^{+}\left|\frac{h_1(0)}{h_2(0)}\right|\right)}{(2|k_0|+1)(\tau-1)\pi}\\
 &\quad+A_{15}+\log \left((2|k_0|+1)\pi\right)+\frac{\log 2}{2}+1= A_{15}\sum\limits_{l=1}^2\left(2T\left(r,h_l\right)+\log \frac{1}{h_l(0)}\right)\\
 &\quad+\frac{2\tau \left(T\left(\tau r,h_2\right)+T\left(\tau r,h_1\right)+\log\left|\frac{1}{h_1(0)}\right|+\log^{+}\left|\frac{h_1(0)}{h_2(0)}\right|\right)}{(2|k_0|+1)(\tau-1)\pi}+A_{15}+\frac{\log 2}{2}+1\\
 &\quad +\log \left((2|k_0|+1)\pi\right)\leq A_{16}\left(T\left(\tau r,h_1\right)+T\left(\tau r,h_2\right)+1\right)\\
   \end{split}
\end{equation}
for $r>r_0.$ Here and in what follows, $A_{16}$ is some large positive constant such that $A_{16}\geq A_{15}.$ This implies that
\begin{equation}\label{eq3.91X22}
m\left(r, \log\frac{h_2}{h_1}\right)\leq  A_{16}\left(T\left(\tau r,h_1\right)+T\left(\tau r,h_2\right)+1\right) \, \, \text{for} \, \, r>r_0.
\end{equation}
 On the other hand, from \eqref{eq3.61X19} and \eqref{eq3.62X19} we have
\begin{equation}\label{eq3.92X23}
T(\tau r, h_l)<A_{6+l} \left(\log^{+}T(B^2_l\tau r,f)+\log^{+}\log^{+}T(B^2_l\tau r,f)+1\right)
\end{equation}
for the positive integer $l$ satisfying $1\leq l\leq 2$ and the positive number $r$ satisfying \eqref{eq3.64X19} and \eqref{eq3.66X19}. Noting that $f$ is a non-constant meromorphic function, we have that $T(r,f)$
is a monotonically increasing, continuous and unbounded function for $0 < r < +\infty$
(cf.\cite[p.8, Corollary 1]{Hayman1964}). Therefore, we substitute \eqref{eq3.92X23} into \eqref{eq3.91X22} for  $l=1$ and $l=2$ respectively, and then we have
\begin{equation}\label{eq3.93X23}
\begin{split}
 &\quad m\left(r, \log\frac{h_2}{h_1}\right)\leq  A_{16}\sum\limits_{l=l}^2A_{6+l} \left(\log^{+}T(B^2_l\tau r,f)+\log^{+}\log^{+}T(B^2_l\tau r,f)+1\right) \\
  &\leq  A_7 A_{16} \left(\log^{+}T(B^2_1\tau r,f)+\log^{+}\log^{+}T(B^2_1\tau r,f)+1\right) \\
    &\quad +  A_8 A_{16} \left(\log^{+}T(B^2_2\tau r,f)+\log^{+}\log^{+}T(B^2_2\tau r,f)+1\right)\\
       &\leq A_{17}\left(\log T(B_3^2\tau r,f)+1\right),
  \end{split}
\end{equation}
for the positive number $r$ satisfying, say,  \eqref{eq3.64X19} and \eqref{eq3.66X19}. Here $B_3$ as in \eqref{eq3.68X19} is the positive constant such that $B_3=\max\{B_1,B_2\}>\tau>1,$ and $A_{17}$ is some large positive constant such that $ A_{17}>\max\{ A_7 A_{16}, A_8 A_{16}\}.$ Finally we substitute \eqref{eq3.68X19}, \eqref{eq3.93X23} and the estimate \eqref{eq2.11X19} for $S(r,f)+S(r,g)$ into \eqref{eq3.56X19}, and then we find a large positive number, say $r_0>1,$ and a large positive constant $A_{18},$ such that for the positive number $r$ satisfying
 \begin{equation}\label{eq3.94X23}
 1<r_0<r<\tau r<B_1\tau r<B^2_1\tau r<B^3_3\tau r<R<+\infty
   \end{equation}
   and
\begin{equation}\label{eq3.95X23}
 1<r_0<r<\tau r<B_2\tau r<B^2_2\tau r<B^3_3\tau r<R<+\infty,
   \end{equation}
we have for $\rho=B^3_3\tau r$ that
\begin{equation}\label{eq3.96X23}
\begin{split}
 \frac {1}{7} T(r,f)&\leq 33\left( A_{12}\log T(B^3_3\tau r,f)+\log r+A_{12}\right)+33A_{17}\left(\log T(B_3^2\tau r,f)+1\right)\\
&\quad +A_4+C_3\left(\log^{+}T(\rho,f)+\log^{+}\log^{+}T(\rho,f)+\log^{+}\frac{\rho}{\rho-r}+1\right)\\
&= 33\left( A_{12}\log T(B^3_3\tau r,f)+\log r+A_{12}\right)+33A_{17}\left(\log T(B_3^2\tau r,f)+1\right) \\
&\quad+A_4+C_3\left(\log T(B^3_3\tau r,f)+\log\log T(B^3_3\tau r,f)+\log \frac{B^3_3\tau r}{B^3_3\tau r-r}+1\right)\\
&\leq A_{18}\left(\log T(B^3_3\tau r,f)+\log r+1\right).\\
\end{split}
\end{equation}
Without loss of generality, we suppose that the positive constant $B_3$ satisfying $B_3=\max\{B_1,B_2\}>\tau>1$ is large enough such that
\begin{equation}\label{eq3.97X23}
\log^{+}\frac{1}{\rho-r}=\log^{+}\frac{1}{B^3_3\tau r-r}=\log ^{+}\frac{1}{r}=0
\end{equation}
for $\rho=B^3_3\tau r$ and the positive number $r$ satisfying \eqref{eq3.94X23} and \eqref{eq3.95X23}. From \eqref{eq3.96X23} and \eqref{eq3.97X23} we have
\begin{equation}\label{eq3.98X23}
\begin{split}
 &\quad  T(r,f)\leq 7A_{18}\left(\log T(B^3_3\tau r,f)+\log r+1\right)\\
&=7A_{18}\left(\log^{+} T(\rho,f)+\log^{+}\frac{1}{\rho-r}+\log^{+}\rho+\log^{+}\frac{1}{r}+1-\log (B^3_3\tau)\right)
\end{split}
\end{equation}
for $\rho=B^3_3\tau r$ and the positive number $r$ satisfying \eqref{eq3.94X23} and \eqref{eq3.95X23}. From \eqref{eq3.98X23} and Lemma \ref{Ku1978} we have
\begin{equation}\label{eq3.99X23}
\begin{split}
T(r,f)&\leq 56A_{18}\left(\log\frac{B^3_3\tau r}{B^3_3\tau r-r}+\log^{+}(2B^3_3\tau r)
+\log^{+}\frac{1}{2B^3_3\tau r}\right)\\
&\quad +20(14A_{18}+1)^2+14A_{18}(1-\log (B^3_3\tau))\leq A_{19}\log r
\end{split}
\end{equation}
for $\rho=B^3_3\tau r,$ $R=2\rho=2B^3_3\tau r$ and the positive number $r$ satisfying \eqref{eq3.94X23} and \eqref{eq3.95X23}. Here $A_{19}$ is a large positive constant such that $A_{19}>56A_{18}.$ From \eqref{eq3.99X23} and Yang-Yi \cite[p.11, Theorem 1.5]{Yang2003} we see that $f$ is reduced to a rational function. This contradicts the known result that  $f$ and $g$ are transcendental meromorphic functions at the beginning of the proof of Theorem \ref{Theorem1.2}. This completes the proof of Theorem \ref{Theorem1.2}.
\vskip 2mm
\par
\section{On a question of Gary G. Gundersen concerning the nonexistence of two distinct and non-constant meromorphic functions sharing three distinct values DM and a fourth value CM }
\vskip 2mm
\par In 1979, G.G.Gundersen proposed the following two questions in Gundersen \cite{Gundersen1979}:
\vskip 2mm
\par
\begin{question}\rm{}(\cite[p.458]{Gundersen1979}) \label{question4.1}
Do there exist two distinct and non-constant meromorphic functions that share two distinct values CM and share the other two distinct values DM, where the four values are four distinct values in the extended complex plane?
\end{question}
\vskip 2mm
\par
\begin{question}\rm{}(\cite[p.458]{Gundersen1979})\label{question4.2}
Do there exist two distinct and non-constant meromorphic functions that share one value CM and share the other three distinct values DM, where the four values are four distinct values in the extended complex plane?
\end{question}
\vskip 2mm
\par In 1983, Gundersen\cite[Theorem 1]{Gundersen1983} proved that if two distinct and non-constant meromorphic functions share two values CM and share two other values IM, where the four values are four distinct values in the extended complex plane, then the two meromorphic functions share all four values CM. This implies that Gundersen\cite[Theorem 1]{Gundersen1983} gave a negative answer to Question \ref{question4.1}. But Question \ref{question4.2} is still open by now. The following question is a special case of Question \ref{question4.2}:
\vskip 2mm
\par
\begin{question}\rm{}\label{question4.3}
Do there exist two distinct and non-constant entire functions that share three distinct finite values DM?
\end{question}
In this direction, we recall the following result in Mues \cite{Mues1989F1} that completely resolved Question \ref{question4.3}:
\vskip 2mm
\par
\begin{theorem}\rm{}(\cite[pp.109-117]{Mues1989F1})\label{TheoremI}
There do not exist two distinct and non-constant entire functions that share three distinct finite values $a_1,$ $a_2,$ $a_3$ DM.
\end{theorem}
\vskip 2mm
\par For the existence of three distinct and non-constant meromorphic functions sharing four distinct values IM, where the four distinct shared values are in the extended complex plane that are neither CM shared values nor DM shared values, we refer to Steinmetz \cite[Theorem 2]{Steinmetz1988}.
For the existence of two distinct and non-constant meromorphic functions sharing four distinct values DM, where the four distinct shared values are in the extended complex plane, we refer to  Reinders \cite[Theorem 1]{Reinders1992}. From Theorem \ref{Theorem1.2} we get the following result that gives a negative answer to Question \ref{question4.2} for the case of meromorphic functions:
\vskip 2mm
\par
\begin{theorem}\rm{}\label{Theorem4.1}
There do not exist two distinct and non-constant meromorphic functions that share $a_1,$ $a_2,$ $a_3$ DM and share $a_4$ CM, where
$a_1,$ $a_2,$ $a_3,$ $a_4$ are four distinct complex values in the extended complex plane.
\end{theorem}
\vskip 2mm
\par
From Example \ref{Example1.1} above and the following example we see that the number of the three distinct DM shared values and the number of one CM shared value in Theorem \ref{Theorem4.1}, in a sense, are best possible:
\vskip 2mm
\par
\begin{example} \rm{}(\cite[p.94]{Steinmetz1988}) Let $\mathcal{P}$ denote the Weierstrass $\mathcal{P}$-function with a pair of primitive periods $2\omega-\omega'$ and $\omega+\omega'.$ Then
\begin{equation}\nonumber
F(z)=\frac{\left(\mathcal{P}(z)-\mathcal{P}(\omega/2)\right)\left(\mathcal{P}(z)-\mathcal{P}(3\omega/2)\right)^2}{\mathcal{P}(z)-\mathcal{P}(\omega)}, \, \, G(z)=F(z+\omega)
\end{equation}
 and $H(z)=F(z+\omega')$ share the value $0,$ $a_1,$ $a_2,$ $\infty$ IM, where $a_1$ and $a_2$ are two distinct complex numbers such that $a_1/a_2$ with $a_1/a_2\neq -1$ is a cube root of unity, while $\omega$ and $\omega'$ are two complex numbers such that $2\omega-\omega'$ and $\omega+\omega'$ are linearly independent on $\Bbb{R}.$ That is to say, for each pair of real numbers $\lambda_1$ and $\lambda_2$ such that $|\lambda_1|+|\lambda_2|\neq 0,$ we have $\lambda_1(2\omega-\omega')+\lambda_2(\omega+\omega')\neq 0.$  Moreover, we can verify that the four values $0,$ $a_1,$ $a_2,$ $\infty$ are neither CM shared values nor DM shared values of $F(z),$ $G(z)$ and $H(z).$
\end{example}
\vskip 2mm
\par
{\bf Declarations.} On behalf of all authors of the present paper, the corresponding author states that there is no conflict of interest.
\vskip 2mm
\par
\def\cprime{$'$}
\providecommand{\bysame}{\leavevmode\hbox to3em{\hrulefill}\thinspace}
\providecommand{\MR}{\relax\ifhmode\unskip\space\fi MR }
\providecommand{\MRhref}[2]{%
\href{http://www.ams.org/mathscinet-getitem?mr=#1}{#2}
}
\providecommand{\href}[2]{#2}

\end{document}